\documentclass[final,english]{siamltex}

\usepackage{babel,amssymb}
\usepackage{graphicx}
\usepackage{amsmath}
\usepackage{bm}
\usepackage[show]{ed}
\usepackage[notcite,notref]{showkeys}
\usepackage{epstopdf}

\newcommand{\punto}{\,\cdot\,}

\newcommand{\smallfrac}[2]{{\textstyle\frac{#1}{#2}}} 
\newcommand{\jump}[1]{[\![#1]\!]}
\newcommand{\ave}[1]{\{\!\!\{#1\}\!\!\}}

\numberwithin{equation}{section}

\title{A Nystr\"om flavored Calder\'on Calculus of order three for two dimensional waves}

\date{\today}
\author{V\'\i ctor Dom\'\i nguez\thanks{Departamento de Ingenier\'\i a
Matem\'atica e Inform\'atica, Universidad P\'ublica de Navarra, 31500 Tudela,
Spain. {\tt victor.dominguez@unavarra.es}. Partially supported by MICINN Project
MTM2010-21037} \and Sijiang L. Lu\thanks{Department of Mathematical Sciences,
University of Delaware, USA. {\tt sjlv@math.udel.edu}} \and Francisco--Javier
Sayas\thanks{Department of Mathematica Sciences, University of Delaware,
Newark, DE 19716, USA. {\tt fjsayas@math.udel.edu}. Partially supported by NSF
grant DMS 1216356.}}

\begin{document}

\maketitle

\begin{abstract}
In this paper we present and test a full discretization of all elements of the Calder\'on 
Calculus (layer potentials and integral operators) for the Helmholtz equation in  smooth 
closed curves in the plane. The resulting integral equations provide approximations of order 
three for all variables involved. Test are shown for a wide array of direct, indirect and 
combined field integral equation at fixed frequency and for a Convolution Quadrature based 
approximation in the time domain.
\end{abstract}


\section{Introduction}

This paper introduces and tests a fully discrete Calder\'on Calculus for two dimensional 
acoustic waves, time-harmonic and transient. We present a novel simultaneous discretization 
of the two layer potentials and four integral operators associated to the Helmholtz equation 
at fixed frequency on a collection of smooth non-intersecting parametrizable closed curves 
in the plane. The method's vocation is simplicity, and we can assert with some confidence, 
that it will not be easy to find another instance of such a simple method of reasonable 
order (order three in all variables, in strong norms), with so little computational and 
programming requirements. We do not make any claims, though, on the ability of this method 
to work on problems at  high frequencies, and we are by no means competitors of 
sophisticated high order methods that look for fine details in complicated geometries. We 
do, however, claim that the set of tools exposed in this paper works for many other integral 
operators (experiments on the Laplace equation have been carried out by the authors as a 
prototyping tool), and we are working on the extension of these ideas to some more general 
problems. It is important to emphasize that we are not discretizing a particular integral 
equation, so we are not worrying on whether one formulation is better than another, or 
whether there are resonances. As we show in the examples, the discrete operators and 
potentials can be used to build any of the best known direct, indirect, and combined field 
integral equations for exterior problems, as well as more complicated systems of integral 
equations for transmission problems. The time domain extension is carried out by using a 
variable complex frequency and the Convolution Quadrature technology of Christian Lubich 
\cite{Lubich:1994, BaSc:2012}.

The method works in a relatively simple way. On a parametric curve, several sets of points 
and normal vectors are sampled. They are harvested using three staggered uniform grids in 
parameter space. One grid is used to create sources and two grids are used for simultaneous 
(averaged) observation. Once each curve is sampled at the discrete level, merging 
information to create a unified discrete set is an easy task. The second step is the 
automatic creation of potentials and operators using direct evaluations of the kernel 
functions: all numerical integration processes are done explicitly in the method, and all 
equations and right-hand sides are fully discrete. The ideas behind these methods go back to 
a very simple quadrature method of order two by Saranen and Schroderus \cite{SaSc:1993}, 
later generalized \cite{CeDoSa:2002} and improved to third order of convergence 
\cite{DoRaSa:2008}. The treatment of the associated hypersingular operator is surprisingly 
simple as well: using the integration by parts formula that is common to Galerkin 
discretizations of the hypersingular integral equation, the paper \cite{DoLuSa:2012a} found 
a collection of fully discrete discretizations of order one and two for this hypersingular 
equation. Some additional work allowed us to put together the first Calder\'on Calculus of 
order two in \cite{DoLuSa:2012}. This set of discrete operators is heavily asymmetric and 
has the disadvantage of requiring sampling of second derivatives of the parametrization of 
the curve due to the evaluation of the double layer operator on its diagonal. The current 
set of methods mixes the discoveries of \cite{DoLuSa:2012} and \cite{DoRaSa:2008} to create 
a discrete set of order three that is even simpler than the order two collection. 

The discrete set is, as a matter of fact, a Nystr\"om (quadrature) discretization of the 
integral operators, avoiding evaluation of singular kernels on the diagonal, looking for 
superconvergent location of observation points, and mixing observation grids to partially 
symmetrize the method, and achieve order three. However, the method can be better understood 
as a full discretization, with carefully chosen low order numerical integration, of a non-
conforming Petrov-Galerkin discretization of the integral operators.  The methods will be 
tested on a wide set of integral equations for exterior and transmission problems, and on a 
time-domain scattering problem. We will also test condition numbers of the different 
formulations and the possibility of using Calder\'on preconditioning.

\paragraph{Some discussion on the literature} Nystr\"om methods
\cite{Nystrom:1930, Atkinson:1997} are the most popular choices for integral equations of 
the second kind.
For integral equations of the second kind with smooth periodic kernels, the trapezoidal rule
gives rise to a very powerful method which converges
superalgebraically \cite[Chapter 12]{Kress:1999}. Periodic weakly singular integral 
equations of the second kind (as those that arise from the Helmholtz equation on smooth 
parametrizable domains in the plane) are also amenable to simple methods with superalgebraic 
or exponential order of convergence \cite[Section 3.5]{CoKr:1998} (see also 
\cite{Kussmaul:1969,Martensen:1963}). A comparison of Nystr\"om methods in the plane has 
been carried out in \cite{HaBaMaYo:2011}. The three dimensional case is much more involved 
and consequently less developed. There are Nystr\"om schemes for equations with weakly 
singular kernels, like those of Bruno and Kunyaski \cite{BrKu:2001,BrDoSa:2013} and Wienert 
\cite{CoKr:1998,Wienert:1990, GrSl:2002}. Finally, the QBX methods  
\cite{KlBaGrON:2012,EpGrKl:2012,HaBaMaYo:2011}, originally designed to compute layer 
potentials close to the boundary, are proving to be useful tools to create Nystr\"om methods 
for weakly singular integral equations in two and three dimensions.

\section{Parametrized Calder\'on Calculus}\label{sec:2}

Let $\mathbf x:\mathbb R \to \Gamma\subset \mathbb R^2$ be a smooth ($\mathbf x\in \mathcal 
C^\infty(\mathbb R)$)  regular ($|\mathbf x'(t)|\neq 0$ for all $t$) $1$-periodic ($\mathbf 
x(t+1)=\mathbf x(t)$ for all  $t$) positively oriented parametrization of a simple ($\mathbf 
x(t)\neq \mathbf x(\tau)$ if $0\le t < \tau <1$) closed curve in the plane. We consider the 
parametrized {\em non-normalized} normal vector field $\mathbf n(t):=(x_2'(t),-x_1'(t))$. 
The curve $\Gamma$ divides the plane into a bounded interior domain $\Omega_-$ and its 
unbounded exterior $\Omega_+$.
Given a function $U:\mathbb R^2\setminus\Gamma\to\mathbb C$ that is smooth enough on both 
sides of the interface $\Gamma$, we will write $U^\pm$ for its restrictions of $\Omega_\pm$. 
Restrictions (traces) and normal derivatives on the boundary will be defined as:
\begin{equation}\label{eq:2.0}
\gamma^\pm U:=U^\pm|_\Gamma \circ\mathbf x, \qquad \partial_{\mathbf n }^\pm U :=\big( (\nabla U^\pm)|_\Gamma\circ\mathbf x\big)\cdot\mathbf n.
\end{equation}
As defined, these restrictions to the boundary define periodic functions and that the normal 
derivative is actually a directional derivative with respect to the non-normalized outward 
pointing normal vector field $\mathbf n$.

Given complex valued sufficiently smooth $1$-periodic functions, we can define the {\em single and double potentials} on $\Gamma$ with the formulas
\begin{subequations}\label{eq:2.1}
\begin{eqnarray}
\big(\mathrm S\, \eta\big)(\mathbf z)&:=& \frac{\imath}{4}\int_{0}^1
H_0^{(1)}(k|\mathbf z-\mathbf x(t)|)\eta(t)\,\mathrm dt,  \\
\big(\mathrm D\, \psi\big)(\mathbf z)&:=& \frac{\imath k}{4}\int_{0}^1
H_1^{(1)}(k|\mathbf z-\mathbf x(t)|)\frac{(\mathbf z-{\bf
x}(t))\cdot\mathbf n(t)}{|\mathbf z-\mathbf x(t)|}\psi(t)\,\mathrm dt.
\end{eqnarray}
\end{subequations}
These functions are defined for all $\mathbf z\in \mathbb R^d\setminus\Gamma$. It is well known (it actually follows from a very simple computation) that for any $\eta, \psi$ the function $U:=\mathrm S \eta+\mathrm D \psi \in \mathcal C^\infty(\mathbb R^2\setminus\Gamma)$ solves
\begin{equation}\label{eq:2.2}
\Delta U+k^2 U=0 \mbox{ in $\mathbb R^2\setminus\Gamma$}, \qquad \lim_{|\mathbf z|\to \infty} |\mathbf z|^{1/2} \Big( \nabla U(\mathbf z)\cdot \big( \smallfrac1{|\mathbf z|}\,\mathbf z\big)-\imath\,k\,U(\mathbf z)\Big) =0,
\end{equation}
that is, $U$ is a radiating solution of the Helmholtz equation. Smoothness of $U$ as we approach the interface $\Gamma$ depends on smoothness of the densities $\eta$ and $\psi$. Reciprocally, given a solution of \eqref{eq:2.2}, we can write
\begin{equation}\label{eq:2.3}
U = \mathrm S \jump{\partial_{\mathbf n} U}-\mathrm D \jump{\gamma U},
\end{equation}
where the jump operators are defined as
\[
\jump{\gamma U}:=\gamma^- U-\gamma^+ U, \qquad \jump{\partial_{\mathbf n} U}:=\partial_{\mathbf n}^-U-\partial_{\mathbf n}^+ U.
\]
Uniqueness of the representation \eqref{eq:2.3} for the solutions of \eqref{eq:2.2} implies the jump relations of potentials
\begin{equation}\label{eq:2.4}
\jump{\gamma \mathrm S\,\eta}=0, \qquad \jump{\partial_{\mathbf n}\mathrm S \,\eta}=\eta, \qquad \jump{\gamma \mathrm D \psi}=-\psi, \qquad \jump{\partial_{\mathbf n}\mathrm D \psi}=0.
\end{equation}
These jump properties motivate the introduction of the four operators on the boundary $\Gamma$:
\begin{subequations}\label{eq:2.5}
\begin{alignat}{4}
&\mathrm V \eta := \ave{\gamma \mathrm S\,\eta}=\gamma^\pm \mathrm S \eta, & \qquad &  \mathrm J\,\eta :=\ave{\partial_{\mathbf n}\mathrm S \,\eta},\\
& \mathrm K\,\psi :=\ave{\gamma\mathrm D\,\psi},  & & \mathrm W\psi:=-\ave{\partial_{\mathbf n}\mathrm D \,\psi}=-\partial_{\mathbf n}^\pm \mathrm D\,\psi,
\end{alignat}
\end{subequations}
where
\[
\ave{\gamma U}:=\smallfrac12 (\gamma^- U+\gamma^+ U), \qquad \ave{\partial_{\mathbf n} U}:=\smallfrac12 (\partial_{\mathbf n}^- U+\partial_{\mathbf n}^+ U).
\]
The operators $\mathrm V$ and $\mathrm K$ are the single and double layer operators respectively,  $\mathrm J$ is   the adjoint double layer operator, and $\mathrm W$ is  the hypersingular operator for the Helmholtz equation. The first three of these operators admit integral expressions:
\begin{subequations}\label{eq:2.6}
\begin{eqnarray} 
({\mathrm  V}\eta)(\tau)&:=&\frac{\imath}{4}\int_0^1H_0^{(1)}(k|\mathbf{x}(\tau)-{
\mathbf  x}(t)|)\eta(t)\,\mathrm dt,\\ 
(\mathrm K\psi)(\tau)&:=&  
\frac{\imath k}{4}\int_0^1 
 H_1^{(1)}(k|\mathbf x(\tau)-\mathbf x(t)|) \frac{(\mathbf x(\tau)-{\bf 
x}(t))\cdot\mathbf n(t)}{|\mathbf x(\tau)-\mathbf x(t)|} 
\psi(t)\,\mathrm dt,\\ 
(\mathrm{J}\eta)(\tau)&:=&  
\frac{\imath k}{4}\int_0^1 
 H_1^{(1)}(k|\mathbf x(\tau)-\mathbf x(t)|) \frac{(\mathbf x(t)-{\bf 
x}(\tau))\cdot\mathbf n(\tau)}{|\mathbf x(\tau)-\mathbf x(t)|} 
\eta(t)\,\mathrm dt.
\end{eqnarray}
It is clear from here that $\mathrm J=\mathrm K^t$. We will keep a different notation though, for reasons that will become apparent when we discretize them in a non--symmetric form. The operator $\mathrm W$ admits an expression in the form of an integro-differential operator:
\begin{equation}\label{eq:2.6d}
\mathrm W\psi:= -( \mathrm V \psi')' -k^2 \mathrm V_{\mathbf n}\psi,
\end{equation}
where
\begin{equation}\label{eq:2.6e}
(\mathrm V_{\mathbf n}\psi)(\tau):=\frac{\imath}{4}\int_0^1H_0^{(1)}(k|\mathbf{x}(\tau)-{\bf 
x}(t)|) \big(\mathbf n(t)\cdot\mathbf n(\tau)\big)\psi(t)\,\mathrm dt.
\end{equation}
\end{subequations}
The representation formula \eqref{eq:2.3} for all radiating solutions of the Helmholtz 
equation \eqref{eq:2.2}, together with the jump properties of the potentials \eqref{eq:2.4} 
and the definitions of  the boundary integral operators by averaging \eqref{eq:2.5}, 
determines a set of rules (a calculus) that generates a diverse collection of 
representation formulas, potential ansatzs, and integral equations associated to the 
solution of interior, exterior and transmission problems for the Helmholtz equation. The 
following matrices of operators
\begin{equation}\label{eq:2.8}
\left[\begin{array}{c} \gamma^\pm \\ \partial_{\mathbf n}^\pm \end{array}\right] \left[\begin{array}{cc} \mathrm D & -\mathrm S \end{array}\right]=
\pm \frac12\left[\begin{array}{cc}\mathrm I & 0 \\ 0 & \mathrm I\end{array}\right]+
\left[\begin{array}{cc} \mathrm K & -\mathrm V \\ -\mathrm W & -\mathrm J\end{array}\right]
\end{equation}
collect the exterior/interior Cauchy values of the layer potentials. They constitute the 
exterior/interior Calder\'on projectors associated to the Helmholtz equation. The 
systematic use of these  potentials and operators to build integral equations will be 
explored in Section \ref{sec:6}.  At this point, let us emphasize the fact that we are 
striving for a full discretization of the entire set of potentials \eqref{eq:2.1} and 
operators \eqref{eq:2.6}, including also discretization of the restriction operators 
\eqref{eq:2.0} that will be needed to sample, at the discrete level, incoming incident 
waves.

\paragraph{The case of multiple scatterers} Assume that $\Gamma_1,\ldots,\Gamma_M$ are pairwise disjoint curves, parametrized as above by smooth $1$-periodic functions $\mathbf x_\ell$. All potentials and operators can be easily defined for vectors of densities $(\eta_1,\ldots,\eta_M)$ and $(\psi_1,\ldots,\psi_M)$. The integral operators then become matrices of integral operators. For instance, we have operators of the form
\[
\frac{\imath}{4}\int_0^1H_0^{(1)}(k|\mathbf x_\ell(\tau)-{
\mathbf  x_m}(t)|)\eta_m(t)\,\mathrm dt, \qquad \ell,m=1,\ldots,M.
\]

\section{Fully discrete method}\label{sec:3}

\subsection{Geometry}\label{sec:3.1}
For one single curve parametrized with $\mathbf x$ as in Section \ref{sec:2}, we proceed as follows. We take a positive integer $N$ and define $h:=1/N$. Next we define the discrete parameter points $t_j:=h\,j$ and the values
\[
\mathbf m_j:=\mathbf x(t_j), \quad \mathbf b_j:=\mathbf x(t_j-h/2), \quad \mathbf n_j:=h\,\mathbf n(t_j), \quad j\in \mathbb Z_N:=\{1,\ldots,N\}.
\]
The notation $\mathbf m_j$ and $\mathbf b_j$ makes reference to midpoints and breakpoints of a boundary element mesh that is implicit to this method (see Section \ref{sec:4}).
In addition to these sampled quantities, we need two index-based functions that provide the next and previous index modulo $N$: the next-index function is $n :\mathbb Z_N \to \mathbb Z_N$ given by
\[
n(j):=\left\{ \begin{array}{ll} j+1, & 1\le j\le N-1,\\
1, & j=N, \end{array}\right.
\]
while $p:=n^{-1}$. Merging geometric information from several curves is easy:
after sampling two curves $\Gamma_1$ and $\Gamma_2$ with $N_1$ and $N_2$ elements 
respectively,  midpoints, breakpoints, and normals are collected in lists with $N=N_1+N_2$ 
elements, by appending the information of $\Gamma_2$ after the information of $\Gamma_1$. 
We then create the next-index and previous-index functions by juxtaposing the two existing 
functions:
\[
n(j)=\left\{\begin{array}{ll} j+1, & 1\le j\le N_1-1,\\ 1, & j=N_1,\\ j+1, & N_1+1\le j\le N_1+N_2-1,\\ N_1+1, & j=N_1+N_2,\end{array}\right.  \qquad p=n^{-1}.
\]
This merging process can be applied to any finite number of curves, each one discretized 
(sampled) with a different number of points. The quantity $h$ appears only at the time of 
collecting information from a particular curve and is incorporated to quantities related to 
first derivatives of the parametrization. However, at the time of merging, $h$ is absent 
from any expression. From this moment on, $n : \mathbb Z_N \to \mathbb Z_N$ is a 
permutation of $\mathbb Z_N$ and $p=n^{-1}.$

\subsection{Discrete potentials}

The discrete version of the single and double layer potentials \eqref{eq:2.1} is defined by using linear combinations of monopoles and dipoles:
\[
\Phi_j(\mathbf z):= \frac\imath4 H^{(1)}_0(k|\mathbf z-\mathbf m_j|) \qquad \mbox{and} \qquad D_j(\mathbf z):=\frac{\imath k}4 H^{(1)}_1 (k|\mathbf z-\mathbf m_j|) \frac{(\mathbf z-\mathbf m_j)\cdot\mathbf n_j}{|\mathbf z-\mathbf m_j|}.
\]
Given two vectors $\boldsymbol\eta=(\eta_1,\ldots,\eta_N)^\top, \boldsymbol\psi =(\psi_1,\ldots,\psi_N)^\top \in \mathbb C^N$, the discrete potentials
\begin{subequations}
\begin{eqnarray}
\mathrm S_h(\mathbf z)\boldsymbol\eta &:=& \sum_{j=1}^N \eta_j \Phi_j(\mathbf z),\\
\label{eq:3.1b}
\mathrm D_h(\mathbf z)\boldsymbol\psi &:=& \smallfrac1{24}\sum_{j=1}^N \psi_j( D_{p(j)}(\mathbf z)+22 D_j(\mathbf z)+D_{n(j)}(\mathbf z))\\
& = & \smallfrac1{24}\sum_{j=1}^N ( \psi_{p(j)}+22 \psi_j+ \psi_{n(j)}) D_j(\mathbf z),
\nonumber
\end{eqnarray}
\end{subequations}
define solutions of \eqref{eq:2.2}.

\paragraph{A quadrature related matrix} Let us consider the $N\times N$ matrix $\mathrm Q$ given by
\begin{equation}\label{eq:3.20}
\mathrm Q_{i,i}=\smallfrac{11}{12}, \qquad \mathrm Q_{i,n(i)}=\mathrm Q_{i,p(i)} = \smallfrac1{24}, \qquad \mathrm Q_{i,j}=0\quad \mbox{otherwise.}
\end{equation}
When the geometry proceeds from a single sampled curve, and therefore the next-index function is just a right-shift modulo $N$, $\mathrm Q$ is the circulant symmetric matrix
\[
\mathrm Q=\frac1{24}\left[\begin{array}{ccccc} 22 & 1 & & & 1 \\ 1 & 22 & 1 & & \\ & \ddots & \ddots & \ddots  \\& &  1 & 22 & 1 \\
1 & & & 1 & 22
\end{array}\right].
\]
In general, $\mathrm Q$ is block diagonal with blocks of the above form, one for each of 
the curves. This matrix is related to a quadrature formula that will be introduced in 
Section \ref{sec:4}. It is clear that \eqref{eq:3.1b} is just a linear combination of 
dipoles, where either the coefficients are premultiplied by the matrix $\mathrm Q$, or the 
dipoles themselves are mixed using this matrix.

\subsection{Observation grids and mixing matrices}\label{sec:3.3}

Since the integral operators in \eqref{eq:2.6} have singularities at $\tau=t$, we are 
forced to use a different discrete set for testing. We start by defining two sets of 
discrete samples. For a single curve parametrized with $\mathbf x$, we use the same $N$  
and $h=1/N$ to define
\[
\mathbf m_j^\pm:=\mathbf x(t_j\pm h/6), \quad \mathbf b_j^\pm:=\mathbf x(t_j-h/2\pm h/6), \quad \mathbf n_j^\pm:=h\, \mathbf x(t_j\pm h/6), \quad j\in \mathbb Z_N.
\]
As in Section \ref{sec:3.1}, observations on finite collections of curves are merged in a simple way. We demand that the number of discretization and observation points on each curve coincides,  although it can be taken to be different on different curves.

Instead of directly averaging values from both possible choices, we will be considering a more general mixture of the two grids. We start with the $N\times N$ matrix $\mathrm P^+=\mathrm P^+(\alpha)$ with elements
\[
\mathrm P^+_{i,i}:=\smallfrac12\alpha, \qquad \mathrm P_{i,p(i)}^+:=\smallfrac12(1-\alpha), \qquad \mathrm P^+_{i,j}=0\quad\mbox{otherwise}.
\]
The parameter $\alpha>0$ will be discussed in Section \ref{sec:5}.
We also let $\mathrm P^-:=(\mathrm P^+)^\top$. For the case of a single curve (when $n$ is the right-shift modulo $N$), we show two particular cases of interest:
\[
\mathrm P^+(\smallfrac56)=\frac1{12}\left[\begin{array}{cccc} 5 & & & 1 \\ 1 & 5 \\ & \ddots & \ddots \\ & & 1 & 5  \end{array}\right], \qquad \mathrm P^+(1)=\smallfrac12\mathrm I.
\]
Given two vectors $\boldsymbol\xi^\pm \in \mathbb  C^N$, it is easy to see that
\begin{equation}\label{eq:3.21}
(\mathrm P^+\boldsymbol\xi^++\mathrm P^-\boldsymbol\xi^-)_i=\smallfrac12 ((1-\alpha)\xi^+_{p(i)} +\alpha \xi_i^- + \alpha\xi_i^++(1-\alpha)\xi_{n(i)}^-).
\end{equation}
Similarly, the $i$-th element of $\mathrm Q(\mathrm P^+\boldsymbol\xi^++\mathrm P^-\boldsymbol\xi^-)=\mathrm P^+\mathrm Q  \boldsymbol\xi^++\mathrm P^-\mathrm Q \boldsymbol\xi^-$, namely
\begin{eqnarray}\label{eq:3.51}
& & \hspace{-2cm}\nonumber
\smallfrac1{48}\Big((1-\alpha) \xi_{p^2(i)}^++\alpha  \xi_{p(i)}^-+ (22-21\alpha) \xi_{p(i)}^++ (21\alpha+1) \xi_i^-\\
& & + (21\alpha+1) \xi_i^++(22-21\alpha) \xi^-_{n(i)}+ \alpha \xi^+_{n(i)}+(1-\alpha)\xi^-_{n^2(i)}\Big),
\end{eqnarray}
is a weighted local average of the values around the index $i$.

The testing part of the discrete Calder\'on Calculus is applied upon an incident wave. At this point, this is just a function $U^{\mathrm{inc}}:\mathbb R^2\to\mathbb C$ that is smooth around the collection of curves, so that we can evaluate
\begin{subequations}\label{eq:3.22}
\begin{eqnarray}
\boldsymbol\beta_0^\pm & :=& - ( U^{\mathrm{inc}}(\mathbf m_1^\pm),\ldots, U^{\mathrm{inc}}(\mathbf m_N^\pm))^\top,\\
\boldsymbol\beta_1^\pm &:=& - (\nabla U^{\mathrm{inc}}(\mathbf m_1^\pm)\cdot\mathbf n_1^\pm, \ldots, \nabla U^{\mathrm{inc}}(\mathbf m_N^\pm)\cdot\mathbf n_N^\pm)^\top.
\end{eqnarray}
We finally define the observation of the incident wave and its normal derivative with
\begin{equation}
\boldsymbol\beta_0:=\mathrm P^+\boldsymbol\beta_0^++\mathrm P^-\boldsymbol\beta^-_0, \qquad \boldsymbol\beta_1:=
\mathrm Q(\mathrm P^+\boldsymbol\beta_1^++\mathrm P^-\boldsymbol\beta^-_1).
\end{equation}
\end{subequations}

\subsection{Discrete operators}

The discrete operators are defined using the geometric elements of Section \ref{sec:3.1} in 
the integration variable and the observation grids of Section \ref{sec:3.3} in the test 
variable. The subscript $h$ will be used to denote discretization. In the case of several 
curves $\{\Gamma_1,\ldots,\Gamma_M\}$, we can consider that $h:=(1/N_1,\ldots,1/N_M)$, although this is not relevant for the exposition of the methods.

Following \eqref{eq:2.6} we define two sets of discrete operators (based on the principal 
sampling of the geometry, tested on both $\pm$ observation grids). We start with the three 
integral operators
\begin{subequations}
\begin{eqnarray}
\mathrm V_{i,j}^\pm &:=& \frac\imath4 H^{(1)}_0(k|\mathbf m_i^\pm-\mathbf m_j|),\\
\mathrm K_{i,j}^\pm &:=& \frac{\imath \,k}{4} H^{(1)}_1(k|\mathbf m_i^\pm-\mathbf m_j|) \frac{(\mathbf m_i^\pm-\mathbf m_j)\cdot\mathbf n_j}{|\mathbf m_i^\pm-\mathbf m_j|},\\
\mathrm J_{i,j}^\pm &:=&\frac{\imath\,k}4 H^{(1)}_1(k|\mathbf m_j-\mathbf m_i^\pm|) \frac{(\mathbf m_j-\mathbf m_i^\pm)\cdot\mathbf n_i^\pm}{|\mathbf m_j-\mathbf m_i^\pm|}.
\end{eqnarray}
Following \eqref{eq:2.6d}, the discretization of $\mathrm W$ separates the discretization of the principal part
\begin{equation}
\widetilde{\mathrm W}_{i,j}^\pm:=\widetilde{\mathrm V}^\pm_{n(i),n(j)}-\widetilde{\mathrm V}^\pm_{n(i),j}-\widetilde{\mathrm V}^\pm_{i,n(j)}+\widetilde{\mathrm V}^\pm_{i,j}, \qquad \widetilde{\mathrm V}^\pm_{i,j}:=\frac\imath4H^{(1)}_0(k|\mathbf b_i^\pm-\mathbf b_j|),
\end{equation}
from the more regular logarithmic term in \eqref{eq:2.6e}
\begin{equation}
\mathrm V_{\mathbf n,i,j}^\pm:=(\mathbf n_i^\pm\cdot\mathbf n_j)\mathrm V_{i,j}^\pm.
\end{equation}
\end{subequations}
If $\mathrm V^\pm_h, \mathrm K^\pm_h, \mathrm J^\pm_h, \widetilde{\mathrm W}^\pm_h,$ and $\mathrm V_{\mathbf n,h}^\pm$ are the above matrices, we define the matrices of the discrete Calder\'on Calculus by
\begin{subequations}
\begin{eqnarray}
\mathrm V_h &:=& \mathrm P^+\mathrm V_h^++\mathrm P^-\mathrm V_h^-,\\
\mathrm K_h &:=& \mathrm P^+\mathrm K^+_h\mathrm Q + \mathrm P^-\mathrm K^-_h\mathrm Q=
(\mathrm P^+\mathrm K^+_h  + \mathrm P^-\mathrm K^-_h)\mathrm Q,\\
\mathrm J_h &:=& \mathrm P^+\mathrm Q \mathrm J_h^++\mathrm P^-\mathrm Q \mathrm J^-_h=
\mathrm Q( \mathrm P^+ \mathrm J_h^++\mathrm P^- \mathrm J^-_h)\\
\mathrm W_h &:=& \mathrm P^+\widetilde{\mathrm W}^+_h+\mathrm P^-\widetilde{\mathrm W}^-_h-k^2 ( \mathrm P^+ \mathrm Q \mathrm V_{\mathbf n,h}^+\mathrm Q +\mathrm P^- \mathrm Q \mathrm V_{\mathbf n,h}^-\mathrm Q )\\
&=&  \mathrm P^+\widetilde{\mathrm W}^+_h+\mathrm P^-\widetilde{\mathrm W}^-_h-k^2\mathrm Q(\mathrm P^+\mathrm V_{\mathbf n,h}^++\mathrm P^-\mathrm V_{\mathbf n,h}^-)\mathrm Q.
\end{eqnarray}
\end{subequations}
The Calder\'on projectors \eqref{eq:2.8} include the action of two identity operators. Both of them will be approximated by the following mass matrix $\mathrm M=\mathrm M(\alpha)$
\begin{equation}\label{eq:3.5}
\mathrm M_{i,i}:=\smallfrac29(1+3\alpha), \quad \mathrm M_{i,p(i)}=\mathrm M_{i,n(i)}:=\smallfrac1{18}(7-6\alpha), \quad \mathrm M_{i,j}=0\quad \mbox{otherwise}.
\end{equation}
The simplest method corresponds to $\alpha=1$. In this case $\mathrm P^\pm=\frac12\,\mathrm 
I$ and, apart from the action of the matrix $\mathrm Q$ (related to quadrature), we are 
just averaging sets of equations on the two grids. However, even in this simple case, the mass matrix has a circulant tridiagonal structure.

\section{From Nystr\"om to Petrov-Galerkin}\label{sec:4}

In this section we reinterpret all the matrices and testing of right-hand sides given in 
Section \ref{sec:3} as non-conforming Petrov-Galerkin method with numerical quadrature. 
This will be done for the case of a single curve, where we are working with a single 
equation and parametric unit interval ($1$-periodic real line). When there are $M$ curves, 
$M$ copies of the unit interval have to be used. The details just became slightly more 
cumbersome, but all the following arguments can be extended readily.

\begin{figure}
\begin{center}
\includegraphics[width=5cm]{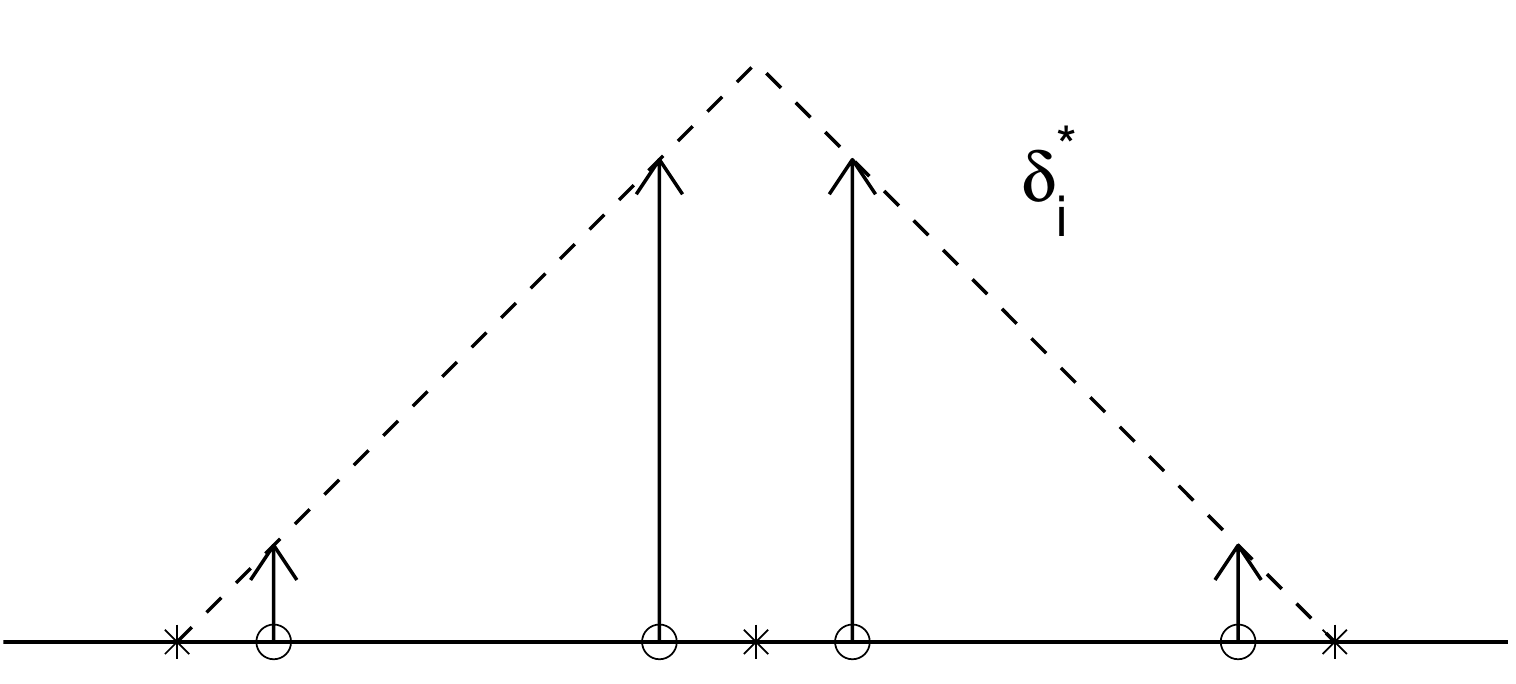}\qquad\includegraphics[width=5cm]{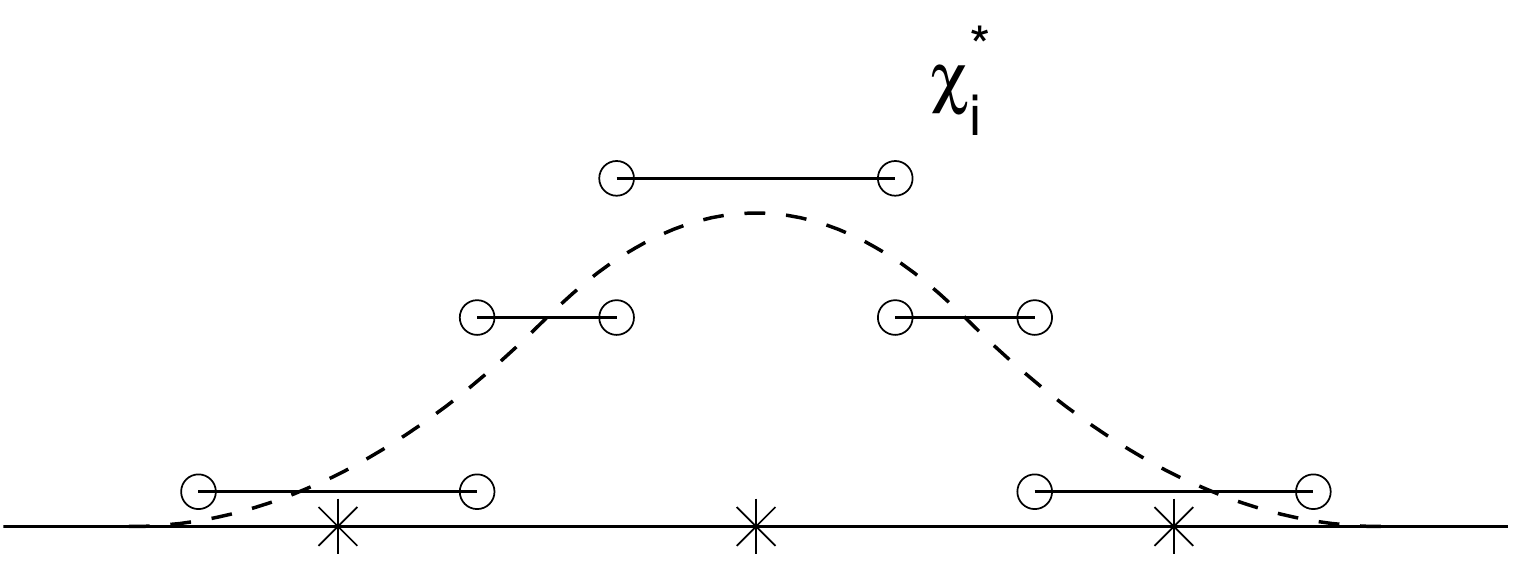}
\end{center}\caption{The shape of the combination of Dirac deltas  $\delta_i^\star$ and the piecewise constant function $\chi_i^\star$. The plot is given for the choice $\alpha=5/6$. A piecewise linear function and a quadratic spline are shown in the background. They are at the origin of the choice of coefficients for the distribution $\delta_i^\star$. }\label{fig:fork}
\end{figure}

\paragraph{Discrete functions and spaces} We start by setting some notation. Given $z\in \mathbb R$, we write $\delta_z$ to denote the $1$-periodic Dirac delta distribution at $z$, that is, the Dirac comb supported on $z+\mathbb Z$. Given an open interval $I$, of length less than one, we write $\chi_I$ to denote the $1$-periodic function that coincides with the characteristic function of $I$ on a unit length interval containing $I$. We then write
\[
\delta_i:=\delta_{t_i}, \quad \delta_i^\pm:=\delta_{t_i\pm h/6}, \quad \chi_i:=\chi_{(t_i-h/2,t_i+h/2)}, \quad \chi_i^\pm :=\chi_{(t_i\pm h/6-h/2, t_i\pm h/6+h/2)}.
\]
Next we define the Dirac fork (see \eqref{eq:3.21} to recognize the corresponding coefficients)
\begin{equation}\label{eq:4.1}
\delta_i^\star:=\smallfrac1{2}\Big((1-\alpha) \delta_{i-1}^++\alpha\delta_i^-+\alpha \delta_i^++(1-\alpha)\delta_{i+1}^-\Big),
\end{equation}
and the ziggurat-shaped piecewise constant functions
\begin{equation}\label{eq:4.2}
\chi_i^\star:=\smallfrac12\Big((1-\alpha) \chi_{i-1}^++\alpha\chi_i^-+\alpha \chi_i^++(1-\alpha)\chi_{i+1}^-\Big).
\end{equation}
Figure \ref{fig:fork} shows the shapes of the basic test functions for the particular case $\alpha=5/6$.
Using momentarily the notation $s_i^\pm:=t_i-h/2\pm h/6$, it is easy to note that, in the sense of periodic distributions,
\[
\smallfrac{\mathrm d}{\mathrm d t}\chi_i^\pm=\delta_{s_i^\pm}-\delta_{s_{i+1}^\pm}
\]
and
\begin{eqnarray*}
\smallfrac{\mathrm d}{\mathrm d t}\chi_i^\star &=& \smallfrac12\Big( (1-\alpha) \delta_{s_{i-1}^+}+\alpha\delta_{s_i^-}+\alpha \delta_{s_i^+}+(1-\alpha)\delta_{s_{i+1}^-}\Big) \\
& & -
 \smallfrac12\Big( (1-\alpha) \delta_{s_{i}^+}+\alpha\delta_{s_{i+1}^-}+\alpha \delta_{s_{i+1}^+}+(1-\alpha)\delta_{s_{i+2}^-}\Big).
\end{eqnarray*}
This shows how, in the same way that characteristic functions arise from integrating two consecutive deltas with opposite signs, the ziggurat functions arise from integrating Dirac forks. Four spaces are relevant for what follows:
\begin{subequations}\label{eq:4.2bis}
\begin{alignat}{4}
& T_h:=\mathrm{span}\{\delta_i\,:\,i\in \mathbb Z_N\}, & \qquad & T_h^\star:=\mathrm{span}\{\delta_i^\star\,:\,i\in \mathbb Z_N\},\\
& S_h:=\mathrm{span}\{\chi_i\,:\,i\in \mathbb Z_N\}, & & S_h^\star:=\mathrm{span}\{\chi_i^\star\,:\,i\in \mathbb Z_N\}.
\end{alignat}
\end{subequations}
Note that $S_h$ is just the space of periodic piecewise constant functions on a uniform mesh with mesh-size $h$ and $\{t_i\}$ as midpoints of the mesh elements.
The $T$ spaces will be non-conforming discretizations of $H^{-1/2}$ Sobolev spaces, while the $S$ spaces are non-conforming approximations of $H^{1/2}$. The $\star$ spaces will do the job of test spaces, while the unscripted spaces will be the trial spaces. 

\paragraph{Interactions of deltas and characteristic functions} We define the actions of deltas with characteristic functions with the formulas
\begin{subequations}\label{eq:4.3}
\begin{alignat}{6}
 \langle \chi_i^\pm,\delta_i\rangle & :=\alpha+\smallfrac{1}{12}&&=: \langle \delta_i^\pm,\chi_i\rangle,\\
\label{eq:4.3b}
\langle \chi_{i+1}^-,\delta_i\rangle=\langle \chi_{i-1}^+,\delta_i\rangle& :=\smallfrac{11}{12}-\alpha&&=:
\langle \delta_{i-1}^+,\chi_i\rangle=\langle \delta_{i+1}^-,\chi_i\rangle,\\
\langle \delta_i,\chi_j^\pm\rangle &:=0&&
=:\langle \delta_j^\pm,\chi_i\rangle, \hspace{2cm}\mbox{otherwise}.
\end{alignat}
\end{subequations}
The {\em otherwise} case above has to be understood modulo $N$. This interaction will be 
explained in Section \ref{sec:5}. It is clear that the first line of \eqref{eq:4.3} 
enforces the second, if we want some kind of consistency of our formulas with respect to 
translations in the origin of the real line. The interactions \eqref{eq:4.3} and the 
definitions \eqref{eq:4.1}, \eqref{eq:4.2} imply that (see \eqref{eq:3.5})
\begin{alignat*}{4}
\langle  \delta_i^\star,\chi_i\rangle=\langle \chi_i^\star,\delta_i\rangle &=\smallfrac29(1+3\alpha)=\mathrm M_{i,i}, \\
 \langle \delta_{i\pm 1}^\star,\chi_i\rangle=\langle \chi_{i\pm 1}^\star,\delta_i\rangle&=\smallfrac1{18}(7-6\alpha)=\mathrm M_{i,i\pm 1}=\mathrm M_{i\pm 1,i}, \\
\langle \delta_i^\star,\chi_j\rangle=\langle \chi_i^\star,\delta_j\rangle &=0, \hspace{1cm}\mbox{otherwise}.
\end{alignat*}
In other words, the matrix $\mathrm M$ is the matrix that represents the `dualities' 
$S_h^\star\times T_h$ and $T_h^\star\times S_h$ if \eqref{eq:4.3} is imposed. It is to be 
noticed that in the simplest case ($\alpha=1$), the interaction of a Dirac delta with a 
characteristic function is forced to be negative on neighboring elements \eqref{eq:4.3b}. 
We will discuss these choices in Section \ref{sec:5}.

\paragraph{First collection of discrete elements} While the angled bracket (linear in both components) is used for the concrete interactions of piecewise constant functions and Dirac delta distributions, from now on we will use curly brackets (linear in both components as well) for the following situations
\[
\{ \delta_z,\phi\}:=\phi(z), \qquad \{\chi_I,\phi\}:=\int_I \phi(t)\mathrm d t.
\]
This can be applied as long as the right-hand side of the expression is meaningful. We can then define the following bilinear forms
\begin{subequations}
\begin{alignat}{6}
T_h^\star\times T_h \ni (\mu_h^\star,\eta_h) & \quad & \longmapsto   \quad &\{ \mu_h^\star,\mathrm V \eta_h\},\\
S_h^\star \times S_h \ni (\phi_h^\star,\psi_h) & \quad & \longmapsto   \quad &  \{ \smallfrac{\mathrm d}{\mathrm dt } \phi_h^\star,\mathrm V \smallfrac{\mathrm d}{\mathrm dt }\psi_h\},
\end{alignat}
as well as the linear map
\begin{equation}
T_h^\star\ni \mu_h^\star \quad \longmapsto \quad \{\mu_h^\star, U^{\mathrm{inc}}\circ\mathbf x\}.
\end{equation}
\end{subequations}
With the given bases for the spaces \eqref{eq:4.2bis}, the bilinear forms produce the matrix $\mathrm V_h$ and $\mathrm P^+\widetilde{\mathrm W}^+_h+\mathrm P^-\widetilde{\mathrm W}^-_h$, while the linear form yields the vector $\boldsymbol\beta_0$.

\paragraph{Look around quadrature} What is missing to get a complete discrete set is the full discretization of the following bilinear forms
\begin{subequations}\label{eq:4.6}
\begin{alignat}{6}
\label{eq:4.6a}
T_h^\star \times S_h \ni (\mu_h^\star,\psi_h) &\quad &\longmapsto \quad & \{ \mu_h^\star,\mathrm K \psi_h\},\\
\label{eq:4.6b}
S_h^\star\times T_h \ni (\phi_h^\star,\eta_h) &\quad &\longmapsto \quad & \{ \phi_h^\star,\mathrm J \eta_h\},\\
\label{eq:4.6c}
S_h^\star\times S_h \ni (\phi_h^\star,\psi_h) &\quad &\longmapsto \quad & \{ \phi_h^\star,\mathrm V_{\mathbf n} \psi_h\},
\end{alignat}
and the linear form
\begin{equation}
\label{eq:4.6d}
S_h^\star\ni\phi_h^\star\quad\longmapsto\quad \{ \phi_h^\star, (\nabla U^{\mathrm{inc}}\circ\mathbf x) \cdot\mathbf n\}.
\end{equation}
\end{subequations}
The elements of the space $S_h^\star$, can be decomposed as sums of elements of the spaces
\[
S_h^\pm:=\mathrm{span}\{\chi_i^\pm\;:\,i\in \mathbb Z_N\}.
\]
Therefore, the practical computation of all elements in \eqref{eq:4.6} can be done if we are able to compute integrals of the form
\begin{equation}\label{eq:4.7}
\int_{a-h/2}^{a+h/2}\mathbf f(t)\cdot\mathbf n(t) \mathrm d t, \qquad \int_{a-h/2}^{a+h/2}\int_{b-h/2}^{b+h/2} m(t,\tau)\,\mathbf n(t)\cdot\mathbf n(\tau)\mathrm dt\mathrm d \tau.
\end{equation}
The approximation of integrals in one variable will be carried out with a three-point formula of order four using points outside the integration interval (see \eqref{eq:3.20})
\begin{equation}\label{eq:4.8}
\int_{a-h/2}^{a+h/2} \phi(t)\mathrm d t\approx \frac{h}{24} (\phi(a-h)+22 \phi(a)+\phi(a+h)).
\end{equation}

\paragraph{Second collection of discrete elements} As already mentioned, the semidiscrete 
elements \eqref{eq:4.6} can be fully discretized once we approximate all integrals of the 
form \eqref{eq:4.7}. For the one variable integrals we use \eqref{eq:4.7} and for the 
double integrals we use the nine-point formula that arises from using \eqref{eq:4.8} in 
each variable. Note that the normal vector appears always in the integration variable and 
that we have defined $\mathbf n_i:=h\, \mathbf n(t_i)$, etc, which means that the value $h$ 
will not appear in any of the resulting expressions. It is then easy to verify that this 
integration process transforms the bilinear forms \eqref{eq:4.6a}--\eqref{eq:4.6c} into 
fully discrete bilinear forms associated to the matrices $\mathrm K_h$, $\mathrm J_h$ and 
$\mathrm Q(\mathrm P^+\mathrm V_{\mathbf n,h}^++\mathrm P^-\mathrm V_{\mathbf 
n,h}^-)\mathrm Q$ respectively. Finally, quadrature on the linear form \eqref{eq:4.6d} leads to the vector $\boldsymbol\beta_1$ in \eqref{eq:3.22}.

\section{Discussion on parameters}\label{sec:5}

There are several choices related to parameters that we next proceed to discuss. The first parameter is the $\pm1/6$ value  that defines the {\em staggered grids} where the spaces $S_h^\pm$ and $T_h^\pm=\mathrm{span}\{\delta_i^\pm:i\in \mathbb Z_N\}$ are defined. These were first discovered in \cite{SaSc:1993} as the optimal choice of the parameter $\varepsilon$ such that the fully discrete method
\[
\sum_{j=1}^N \log |\mathbf x(t_i-\varepsilon/h)\!-\!\mathbf x(t_j)|\,\lambda_j = g(t_i-\varepsilon h) \qquad i=1,\ldots,N
\]
provides a second order approximation of the parametrized Symm's equation
\[
\int_0^1 \log|\mathbf x(\tau)-\mathbf x(t)|\,\lambda(t)\mathrm d t=g(\tau).
\]
All other choices yield methods of order one, except $\varepsilon=0$ which is not 
practicable and $\varepsilon=1/2$ which gives an unstable method. (Note that 
$\varepsilon+1$ leads to the same method as $\varepsilon$.)
With different techniques, these optimal choices were rediscovered in \cite{CeDoSa:2002}, 
where the method was shown to work for more complicated logarithmic kernels (such as the 
one for the Helmholtz equation), and where it was shown that $\varepsilon=\pm1/6$ were the 
only two values that led to second order methods. In fact, after some simplification, the 
leading term of the expansion in \cite[Proposition 16]{CeDoSa:2002} is formally the 
quadrature error (the expansion holds in some Sobolev norm)
\[
\int_0^1\!\! \log_{\#}(t\!-\!\tau)\, u(\tau)\,\mathrm d \tau\!-\!h\sum_{j=1}^N \log_{\#}(t\!-\!t_j) u(t_j)= h\, C\, (\log 4+\log_{\#}(t/h))\, u(t) + \mathcal O(h^2)
\]
in terms of the periodic logarithmic function $\log_{\#}(t):=\log(\sin^2(\pi t)).$
(Note that this was also studied in \cite{CeSa:1999}, where the $\log_\#$ in the first 
order coefficient was not identified, although its graph was given.) This shows clearly 
that the best observation points for this quadrature error are those canceling the order 
one coefficient, namely, the points $t=ih\pm \frac16 h$, which are exactly the points that 
are used in our fully discrete methods.
Only very recently \cite{DoLuSa:2012a}, it was discovered (by the authors of the current 
paper), that the same structure could be used to find a Nystr\"om discretization of the 
hypersingular operator written in integrodifferential form. In its turn, this led to the 
construction of two fully discrete Calder\'on Calculus of order two (one for each of $\varepsilon=\pm1/6$) in \cite{DoLuSa:2012}.

The values $(\frac{11}{12},\frac{1}{24})$ of the matrix $\mathrm Q$ come from the {\em 
look-around quadrature} formula \eqref{eq:4.8}. The need for using points around the 
integration interval in quadratures is related to asymptotic behavior of the discretization 
errors: we want to have quadrature of sufficiently high order, but we do not want to 
introduce any more relative distances between points, since they would trigger first order 
asymptotic errors through the function $C(\varepsilon):=\log 4+\log_{\#}(\varepsilon)$. We 
believe that the formula \eqref{eq:4.8} might be new, but it has to be said that it has 
been derived in the same spirit as formulas in \cite{CrSa:1998} and \cite{HsKoWe:1980}, 
trying to keep fixed relative distances between integration points at the price of using 
points outside the integration interval.

The following set of parameters is given by the definition of the {\em fork} \eqref{eq:4.1} 
and the {\em ziggurat} \eqref{eq:4.2}, that is, they correspond to the matrices $\mathrm 
P^\pm(\alpha)$. The choice $\alpha=\frac56$,  was first discovered in \cite{DoRaSa:2008}, 
applied just to the single layer operator $\mathrm V$. The choice of parameters is 
motivated by the figure of the Dirac deltas fitting in a triangular shape (a basis function 
for the space of continuous piecewise linear functions). This is due to the origin of the 
method based on a variant of the qualocation methods of Ian Sloan \cite{Sloan:2000}. In 
particular, the stability analysis for the corresponding matrix $\mathrm V_h$ (in form of 
an inf-sup condition \cite[Proposition 10]{DoRaSa:2008}) is essentially outsourced to the 
work of Chandler and Sloan on qualocation methods \cite{ChSl:1990}. 
A nice feature of the particular Dirac fork $\alpha=\frac56$, following the shape of a hat 
function, is that its antiderivatives have the shape of the ziggurat, which mimics that 
shape of a B-spline of degree two, as corresponds to antiderivatives of hat functions.

The {\em interactions of Dirac deltas and characteristic functions} can be expressed either with simple elements
\begin{subequations}\label{eq:5.2}
\begin{alignat}{6}
& \langle \chi_i^\pm,\delta_i\rangle:=\gamma, & \quad & \langle \chi_{i+1}^-,\delta_i\rangle=\langle \chi_{i-1}^+,\delta_i\rangle:=1-\gamma, & \quad & \langle \delta_i,\chi_j^\pm\rangle:=0,\quad\mbox{otherwise},\\
& \langle \delta_i^\pm,\chi_i\rangle:=\gamma, & & \langle \delta_{i-1}^+,\chi_i\rangle=\langle \delta_{i+1}^-,\chi_i\rangle:=1-\gamma, & & \langle \delta_j^\pm,\chi_i\rangle:=0, \quad\mbox{otherwise},
\end{alignat}
\end{subequations}
or with the composite actions of forks over simple characteristics and ziggurats over simple deltas (that is, with the elements of the {\em mass matrix} $\mathrm M$):
\begin{subequations}\label{eq:5.3}
\begin{alignat}{4}
\langle  \delta_i^\star,\chi_i\rangle=\langle \chi_i^\star,\delta_i\rangle& =1-2\rho, \\
 \langle \delta_{i\pm 1}^\star,\chi_i\rangle=\langle \chi_{i\pm 1}^\star,\delta_i\rangle&=\rho,\\
\langle \delta_i^\star,\chi_j\rangle=\langle\chi_i^\star,\delta_j\rangle& =0,\quad& \mbox{otherwise}.
\end{alignat}
\end{subequations}
It is clear that, given the parameter $\alpha$ in \eqref{eq:4.1}-\eqref{eq:4.2}, $\gamma$ 
determines $\rho$ and vice versa. What is less obvious, and we will try to explain next, is 
that $\alpha$ (and the choice of the quadrature rule), actually determines both sets of 
coefficients: $\rho=\frac{1}{18}(7-6\alpha) $ and $\gamma=\frac{1}{12}(1+12\alpha).$
We start this argument with a simple computation. Let $T_c$ be the translation operator $T_c \lambda:=\lambda(\cdot-c)$ and consider two collections of formal finite difference operators
\[
\underline\Delta_h^\alpha:=\frac{1-\alpha}2(T_{-\frac56 h}+T_{\frac56h})+\frac\alpha2 (T_{-\frac16h}+T_{\frac16h}), \quad \Delta_h^\beta:=\beta(T_{-h}+T_h)+(1-2\beta) T_0.
\]
With this notation we can write $\delta_i^\star=\underline\Delta_h^\alpha\delta_i$, $\chi_i^\star=\underline\Delta_h^\alpha\chi_i$ and \eqref{eq:4.8} becomes
\begin{eqnarray}
\nonumber
\{\chi_{(a-h/2,a+h/2)},\phi\}=\int_{a-h/2}^{a+h/2}\phi(t)\mathrm d t & \approx & \smallfrac{h}{24}(\phi(a-h)+22\phi(a)+\phi(a+h))\\
&=&
h\{\delta_a,\Delta_h^{1/24} \phi\}=
h\{ \Delta_h^{1/24}\delta_a,\phi\}.
\label{eq:5.4}
\end{eqnarray}
Similarly, the action of the composite difference operator $\underline\Delta_h^\alpha \Delta_h^{1/24}=\Delta_h^{1/24}\underline\Delta_h^\alpha$ is given by the expression
\[
\smallfrac1{48} \Big( (21\alpha+1) (T_{-\frac16h}+T_{\frac16h})+(22-21\alpha)(T_{-\frac56h}+T_{\frac56h})+\alpha (T_{-\frac76h}+T_{\frac76h})+(1-\alpha)(T_{-\frac{11}6h}+T_{\frac{11}6h})\Big)
\]
(recall \eqref{eq:3.51}). A simple computation then shows that 
\begin{equation}\label{eq:5.5}
\Delta_h^{1/24} \underline\Delta_h^\alpha \eta-\Delta_h^\rho \eta = \mathcal O(h^4) \smallfrac{\mathrm d^4}{\mathrm d t^4}\quad\Longleftrightarrow\quad \rho=\smallfrac{1}{18}(7-6\alpha).
\end{equation}
Let us now try to justify why \eqref{eq:5.5} is relevant. Imagine that we want to solve the trivial equation
$\lambda =\partial_{\mathbf n} U^{\mathrm{inc}}$
with our class of methods. The non-conforming Petrov-Galerkin approximation of this equation is
\begin{equation}\label{eq:5.7}
\lambda_h \in T_h, \qquad \langle \chi_i^\star,\lambda_h\rangle=\rho (\lambda_{i-1}+\lambda_{i+1})+(1-2\rho) \lambda_i=\{\chi_i^\star,\partial_{\mathbf n} U^{\mathrm{inc}}\} \quad \forall i.
\end{equation}
The fully discrete method consists of separating $\chi_i^\star$ into its $\pm$ parts and then using quadrature on each side. This leads to the following argument (see \eqref{eq:5.4}):
\[
\{\chi_i^\star,\partial_{\mathbf n} U^{\mathrm{inc}}\} = \{ \underline\Delta_h^\alpha\chi_i,\partial_{\mathbf n} U^{\mathrm{inc}}\}=\{
\chi_i,\underline\Delta_h^\alpha\partial_{\mathbf n} U^{\mathrm{inc}}\}\approx  h \{\delta_i, \Delta_h^{1/24}\underline\Delta_h^\alpha\partial_{\mathbf n} U^{\mathrm{inc}}\}.
\]
The fully discrete realization of $\lambda =\partial_{\mathbf n} U^{\mathrm{inc}}$ is then given by
\begin{equation}\label{eq:5.8}
\lambda_h=\sum_j \lambda_j \delta_j, \qquad \rho (\lambda_{i-1}+\lambda_{i+1})+(1-2\rho) \lambda_i=h\{ \delta_i, \Delta_h^{1/24}\underline\Delta_h^\alpha\partial_{\mathbf n} U^{\mathrm{inc}}\}, \quad \forall i.
\end{equation}
A dimensional look at \eqref{eq:5.8} shows how the unknowns $\lambda_j$ are trying to 
approximate $h \lambda(t_j)$. The consistency error for equations \eqref{eq:5.8} is then 
obtained when plugging in $h\lambda=h \partial_{\mathbf n}U^{\mathrm{inc}}$ in the left 
hand side of the discrete equations and subtracting the right-hand side:
$h ( (\Delta_h^\rho\lambda)(t_i)-(\Delta_h^{1/24}\underline\Delta_h^\alpha\lambda)(t_i)). $
This takes us back to \eqref{eq:5.5}.

\section{Building equations using the discrete calculus}\label{sec:6}

We show here how to write integral equations for boundary value problems associated to the exterior Helmholtz equation:
\[
\Delta U+k^2 U=0 \quad \mbox{in $\Omega_+$}, \qquad \partial_r U-\imath\,k\,U=o(r^{-1/2})\mbox{ at infinity}.
\]
All formulations will be given directly at the discrete level. Here $\Omega_+$ is the exterior of a collection of smooth closed curves with non-intersecting interiors.

\paragraph{Dirichlet problem} With a boundary condition
$\gamma U+\gamma U^{\mathrm{inc}}=0$, 
we can try four different formulations. In all cases, the trace of the incident wave is tested using \eqref{eq:3.22}.  In the indirect formulations we have to give the integral equation and the potential representation. A single layer potential leads to an integral equation of the first kind
\begin{equation}\label{eq:6.1}
\mathrm V_h\boldsymbol\eta=\boldsymbol\beta_0 \qquad \mbox{and}\qquad U_h=\mathrm S_h(\punto)\boldsymbol\eta,
\end{equation}
while a double layer potential leads to an integral equation of the second kind
\begin{equation}\label{eq:6.2}
\smallfrac12\mathrm M\boldsymbol\psi+\mathrm K_h \boldsymbol\psi=\boldsymbol\beta_0 \qquad \mbox{and} \qquad U_h=\mathrm D_h(\punto)\boldsymbol\psi.
\end{equation}
In the direct formulations, we have a representation formula in terms of discrete Cauchy data:
\begin{equation}\label{eq:6.3}
U_h=\mathrm S_h(\cdot)\boldsymbol\lambda-\mathrm D_h(\cdot)\boldsymbol\varphi.
\end{equation}
Here  $\boldsymbol\lambda$ can be found using one of two integral equations and $\boldsymbol\varphi$ will be derived by projecting data. We can use an integral equation of the first kind
\begin{equation}\label{eq:6.4}
\mathrm V_h\boldsymbol\lambda=-\smallfrac12\mathrm M\boldsymbol\varphi+\mathrm 
K_h\boldsymbol\varphi, \qquad \mbox{where} \qquad \mathrm M\boldsymbol\varphi=\boldsymbol\beta_0,
\end{equation}
or an equation of the second kind
\begin{equation}\label{eq:6.5}
\smallfrac12\mathrm M\boldsymbol\lambda+\mathrm J_h\boldsymbol\lambda=-\mathrm W_h\boldsymbol\varphi , \qquad \mbox{where} \qquad \mathrm M\boldsymbol\varphi=\boldsymbol\beta_0.
\end{equation}
In both cases, $\lambda_i \approx \nabla U(\mathbf m_i)\cdot\mathbf n_i.$
\paragraph{Neumann problem} Consider now a boundary condition $\partial_{\mathbf n} U+\partial_{\mathbf n} U^{\mathrm{inc}}=0,$
and test the incident wave as in \eqref{eq:3.22} to produce a vector $\boldsymbol\beta_1$. There are two possible indirect formulations: with the single layer potential
\begin{equation}\label{eq:6.6}
-\smallfrac12\mathrm M\boldsymbol\eta+\mathrm J_h\boldsymbol\eta=\boldsymbol\beta_1 \qquad \mbox{and}\qquad U_h=\mathrm S_h(\punto)\boldsymbol\eta
\end{equation}
and with the double layer potential
\begin{equation}\label{eq:6.7}
\mathrm W_h\boldsymbol\psi=-\boldsymbol\beta_1 \qquad \mbox{and}\qquad U_h=\mathrm D_h(\punto)\boldsymbol\psi.
\end{equation}
The direct formulations use the representation formula \eqref{eq:6.3} and either the equations
\begin{equation}\label{eq:6.8}
-\smallfrac12\mathrm M\boldsymbol\varphi+\mathrm K_h\boldsymbol\varphi =\mathrm V_h\boldsymbol\lambda, \qquad \mbox{where}\qquad \mathrm M\boldsymbol\lambda=\boldsymbol\beta_1,
\end{equation}
or
\begin{equation}\label{eq:6.9}
\mathrm W_h\boldsymbol\varphi=-\smallfrac12\mathrm M\boldsymbol\lambda-\mathrm J_h\boldsymbol\lambda, \qquad \mbox{where}\qquad \mathrm M\boldsymbol\lambda=\boldsymbol\beta_1.
\end{equation}
In the direct representation $\varphi_i\approx \gamma U(\mathbf m_i)$.

\paragraph{Combined potentials} If $-k^2$ is a Dirichlet eigenvalue of the Laplace operator 
in the interior domain $\Omega_-$, then equations \eqref{eq:6.1}, \eqref{eq:6.4}, 
\eqref{eq:6.6} and \eqref{eq:6.8} are approximations of not uniquely solvable problems. 
Similarly, if $-k^2$ is a Neumann eigenvalue, all other four equations break down. Well 
posed equations for all frequencies can be found using a combined field integral 
representation:
\begin{equation}\label{eq:6.10}
U_h=(\mathrm D_h(\punto)-\imath\,k\,\mathrm S_h(\punto))\boldsymbol\eta,
\end{equation}
leading to
\begin{equation}\label{eq:6.11}
\smallfrac12\mathrm M\boldsymbol\eta+\mathrm K_h\boldsymbol\eta-\imath\,k\,\mathrm V_h\boldsymbol\eta=\boldsymbol\beta_0
\end{equation}
for the Dirichlet problem, and
\begin{equation}\label{eq:6.12}
-\mathrm W_h\boldsymbol\eta+\imath\,k\,\smallfrac12\mathrm M\,\boldsymbol\eta-\imath\,k\,\mathrm J_h\boldsymbol\eta=\boldsymbol\beta_1
\end{equation}
for the Neumann problem. Direct formulations based on combined field equations can also be derived using the arguments of the Burton-Miller integral equation. 

\section{Experiments in the frequency domain}
Let $\Gamma_1$ be parametrized by
\begin{equation}\label{eq:7.1}
t\mapsto (\smallfrac1{10},\smallfrac2{10})+\smallfrac1{\sqrt2} ((1+\cos^2(2\pi t))\cos(2\pi t),(1+\sin^2(2\pi t))\sin(2\pi t))\left( \begin{array}{cc} 1 & -1 \\ 1 &1 \end{array}\right),
\end{equation}
and let $\Gamma_2$ be the ellipse parametrized by
$t\mapsto (4,5)+(\cos(2\pi t), 2\sin(2\pi t)).$
Discretization will be led by a single parameter $N$: we will take $2N$ points on $\Gamma_1$ and $N$ points on $\Gamma_2$. We fix the wave number  $k=3$ and consider a source point solution
\begin{equation}\label{eq:7.2}
U(\mathbf x)=\frac\imath 4 H^{(1)}_0(k|\mathbf x-\mathbf x_0|)\qquad\mbox{with}\qquad  \mathbf x_0:=(\smallfrac1{10},\smallfrac2{10}).
\end{equation}
Since the point $\mathbf x_0$ is in the interior of $\Gamma_1$, using $U^{\mathrm{inc}}=-U$ 
as incident wave, will give $U$ as exact solution of the corresponding exterior problem. 
The boundaries of the scatterers are thus acting as transparent screens.
We will measure errors
\begin{equation}\label{eq:7.3}
\mathrm E_N^{\mathrm{ext}}:=\max_{\mathbf z\in \mathrm{Obs}}| U(\mathbf z)-U_h(\mathbf z)|, \qquad \mathrm{Obs}=\{(0,4),(4,0),(-4,2),(2,-4)\}.
\end{equation}
For direct methods involving the computation of $\boldsymbol\lambda$, we will compute
\[
\mathrm E_N^\lambda:=N\, \max_j | \lambda_j-\nabla U(\mathbf m_j)\cdot\mathbf n_j|.
\]
The rescaling factor $N$ is due to the fact that $|\mathbf n_j|$ is proportional to $h$, instead of being of order one. For direct methods involving  $\boldsymbol\varphi$, we will compute
\[
\mathrm E_N^\varphi:=\max_j | \phi_j-U(\mathbf m_j)|, \quad \mbox{where}\quad\boldsymbol\phi=\mathrm Q \boldsymbol\varphi.
\]
Note that the effective approximation of the trace in the discrete potential 
\eqref{eq:3.1b} is not $\boldsymbol\varphi$ but $\boldsymbol\phi=\mathrm Q 
\boldsymbol\varphi$, which justifies our choice for the latter to compute norms of errors. 
It is clear that $\mathrm E_N^{\mathrm{ext}}$ measures the error of a smoothing postprocess 
and, as such, will benefit from weak superconvergence properties. On the other hand, the 
errors for the quantities on the boundary are measured in uniform norm. We will show that 
in all the experiments and for all the quantities, the errors are $\mathcal O(N^{-3})$. 
Experimental orders of convergence are computed using errors on two consecutive meshes.

\paragraph{First round of experiments} We first test all the formulations of Section 
\ref{sec:6} using the above geometry and exact solution. In all of them we test the 
simplest method ($\alpha=1$) and the method that generalizes the fork distribution in 
\cite{DoRaSa:2008} ($\alpha=5/6$), for which there is partial theoretical justification. 
Tables \ref{table:eq1} to \ref{table:eq10}  show convergence of order three in all 
measurable errors. Note that the method for $\alpha=5/6$ is almost invariably slightly 
better than the method for $\alpha=1$. The errors are displayed in Tables \ref{table:eq1} 
to \ref{table:eq10}, corresponding to the ten integral equations given in Section \ref{sec:6}.

\begin{table}[htb]
\begin{center}
\begin{tabular}{rcccc}\hline
 $N$&$\alpha = 5/6$&e.c.r&$\alpha=1$&e.c.r\\[0.5ex] \hline
10&2.0504$E(-001)$&           &2.1788$E(-001)$ &\\
20&4.2900$E(-003)$&5.5788&6.8665$E(-003)$ &4.9879\\
40&4.2678$E(-004)$&3.3294&7.6927$E(-004)$ &3.1580\\
80&5.0466$E(-005)$&3.0801&9.3497$E(-005)$ &3.0405\\
160&6.2217$E(-006)$&3.0199&1.1603$E(-005)$ &3.0104\\
320&7.7503$E(-007)$&3.0050&1.4477$E(-006)$ &3.0027\\
640&9.6795$E(-008)$&3.0012&1.8087$E(-007)$ &3.0007\\
\hline
\end{tabular}
\end{center}
\caption{Errors $\mathrm E_N^{\mathrm{ext}}$ for equation \eqref{eq:6.1} (indirect, single layer, Dirichlet). }\label{table:eq1}
\end{table}

\begin{table}[htb]
\begin{center}
\begin{tabular}{rcccc}\hline
 $N$&$\alpha = 5/6$&e.c.r&$\alpha=1$&e.c.r\\[0.5ex] \hline
10&1.0885$E(-001)$&           &1.1905$E(-001)$ &\\
20&2.1132$E(-004)$&9.0086&6.1488$E(-004)$ &7.5971\\
40&1.4713$E(-005)$&3.8443&4.0943$E(-005)$ &3.9086\\
80&1.5695$E(-006)$&3.2288&3.1627$E(-006)$ &3.6944\\
160&1.8971$E(-007)$&3.0484&2.8519$E(-007)$ &3.4712\\
320&2.3782$E(-008)$&2.9959&2.9196$E(-008)$ &3.2881\\
640&2.9942$E(-009)$&2.9896&3.2775$E(-009)$ &3.1551\\
\hline
\end{tabular}
\end{center}
\caption{Errors $\mathrm E_N^{\mathrm{ext}}$ for equation \eqref{eq:6.2} (indirect, double layer, Dirichlet). }\label{table:eq2}
\end{table}

\begin{table}[htb]
\begin{center}
\begin{tabular}{rcccc}\hline
 $N$&$\alpha = 5/6$&e.c.r&$\alpha=1$&e.c.r\\[0.5ex] \hline
10&1.1492$E(-001)$&           &1.2300$E(-001)$ &\\
20&9.5390$E(-004)$&6.9125&1.1919$E(-003)$ &6.6892\\
40&1.1902$E(-004)$&3.0026&1.3916$E(-004)$ &3.0985\\
80&1.4778$E(-005)$&3.0097&1.7282$E(-005)$ &3.0095\\
160&1.8395$E(-006)$&3.0060&2.1610$E(-006)$ &2.9995\\
320&2.2948$E(-007)$&3.0029&2.7039$E(-007)$ &2.9986\\
640&2.8657$E(-008)$&3.0014&3,3822$E(-008)$ &2.9990\\
\end{tabular}
\begin{tabular}{rcccc}\hline
 $N$&$\alpha = 5/6$&e.c.r&$\alpha=1$&e.c.r\\[0.5ex] \hline
10&4.5613$E(+000)$&           &4.6869$E(+000)$ &\\
20&2.3802$E(-001)$&4.2603&3.9297$E(-001)$ &3.5761\\
40&1.9732$E(-002)$&3.5925&4.3200$E(-002)$ &3.1853\\
80&2.3458$E(-003)$&3.0724&5.3704$E(-003)$ &3.0079\\
160&2.8639$E(-004)$&3.0340&6.6578$E(-004)$ &3.0119\\
320&3.5581$E(-005)$&3.0088&8.3179$E(-005)$ &3.0007\\
640&4.4405$E(-006)$&3.0023&1.0395$E(-005)$ &3.0003\\
\hline
\end{tabular}
\end{center}
\caption{Errors $\mathrm E_N^{\mathrm{ext}}$ and $\mathrm E_N^\lambda$ for equation 
\eqref{eq:6.4} with exterior solution computed using \eqref{eq:6.3} (direct, weakly 
singular integral equation, Dirichlet). The upper table corresponds to $\mathrm 
E_N^{\mathrm{ext}}$ and the lower table corresponds to $\mathrm E_N^\lambda$.}\label{table:eq3}
\end{table}

\begin{table}[htb]
\begin{center}
\begin{tabular}{rcccc}\hline
 $N$&$\alpha = 5/6$&e.c.r&$\alpha=1$&e.c.r\\[0.5ex] \hline
10&3.4080$E(-001)$&           &3.6686$E(-001)$ &\\
20&1,9862$E(-002)$&4.1008&2.1362$E(-002)$ &4.1021\\
40&1.2691$E(-003)$&3.9682&1.4680$E(-003)$ &3.8631\\
80&8.3324$E(-005)$&3.9289&1.0955$E(-004)$ &3.7441\\
160&6.1749$E(-006)$&3.7542&1.2007$E(-005)$ &3.1896\\
320&5.7108$E(-007)$&3.4347&1.4394$E(-006)$ &3.0603\\
640&6.7160$E(-008)$&3.0880&1.7631$E(-007)$ &3.0293\\
\end{tabular}
\begin{tabular}{rcccc}\hline
 $N$&$\alpha = 5/6$&e.c.r&$\alpha=1$&e.c.r\\[0.5ex] \hline
10&2.5186$E(+001)$&           &2.8815$E(+001)$ &\\
20&1.5341$E(+000)$&4.0371&1.7131$E(+000)$ &4.0722\\
40&1.2257$E(-001)$&3.6457&1.5322$E(-001)$ &3.4830\\
80&9.2022$E(-003)$&3.7355&1.3286$E(-002)$ &3.5276\\
160&7.9479$E(-004)$&3.5355&1.4155$E(-003)$ &3.2305\\
320&7.9863$E(-005)$&3.3150&1.6693$E(-004)$ &3.0840\\
640&9.3918$E(-006)$&3.0880&2.0384$E(-005)$ &3.0338\\
\hline
\end{tabular}
\end{center}
\caption{Errors $\mathrm E_N^{\mathrm{ext}}$ and $\mathrm E_N^\lambda$ for equation 
\eqref{eq:6.5} with exterior solution computed with \eqref{eq:6.3} (direct, second kind 
integral equation, Dirichlet). The upper table corresponds to $\mathrm E_N^{\mathrm{ext}}$ 
and the lower table corresponds to $\mathrm E_N^\lambda$.}\label{table:eq4}
\end{table}

\begin{table}[htb]
\begin{center}
\begin{tabular}{rcccc}\hline
 $N$&$\alpha = 5/6$&e.c.r&$\alpha=1$&e.c.r\\[0.5ex] \hline
10&2.0581$E(-001)$&           &2.2349$E(-001)$ &\\
20&8.9154$E(-005)$&1.1173&2.6185$E(-004)$ &9.7372\\
40&1.3500$E(-005)$&2.7233&1.7638$E(-005)$ &3.8920\\
80&1.4812$E(-006)$&3.1881&1.3681$E(-006)$ &3.6885\\
160&1.6462$E(-007)$&3.1696&1.3916$E(-007)$ &3.2974\\
320&1.9224$E(-008)$&3.0981&1.6450$E(-008)$ &3.0806\\
640&2.3189$E(-009)$&3.0514&2.1447$E(-009)$ &2.9392\\
\hline
\end{tabular}
\end{center}
\caption{Errors $\mathrm E_N^{\mathrm{ext}}$ for equation \eqref{eq:6.6} (indirect, single layer, Neumann). }\label{table:eq5}
\end{table}

\begin{table}[htb]
\begin{center}
\begin{tabular}{rcccc}\hline
 $N$&$\alpha = 5/6$&e.c.r&$\alpha=1$&e.c.r\\[0.5ex] \hline
10&1.4507$E(-001)$&           &1.3385$E(-001)$ &\\
20&1.8995$E(-002)$&2.9330&1.9220$E(-002)$ &2.7999\\
40&9.3066$E(-004)$&4.3512&9.3855$E(-004)$ &4.3560\\
80&6.1122$E(-005)$&3.9285&6.1330$E(-005)$ &3.9358\\
160&4.3175$E(-006)$&3.8234&4.3356$E(-006)$ &3.8223\\
320&3.3804$E(-007)$&3.6749&4.2660$E(-007)$ &3.3453\\
640&3.0335$E(-008)$&3.4782&5.1771$E(-008)$ &3.0427\\
\hline
\end{tabular}
\end{center}
\caption{Errors $\mathrm E_{N}^{\mathrm{ext}}$ for equation \eqref{eq:6.7} (indirect, double layer, Neumann). }\label{table:eq6}
\end{table}

\begin{table}[htb]
\begin{center}
\begin{tabular}{rcccc}\hline
 $N$&$\alpha = 5/6$&e.c.r&$\alpha=1$&e.c.r\\[0.5ex] \hline
10&1.8360$E(-001)$&           &2.1658$E(-001)$ &\\
20&3.2147$E(-003)$&5.8357&5.4219$E(-003)$ &5.3199\\
40&3.2038$E(-004)$&3.3268&5.9516$E(-004)$ &3.1874\\
80&3.7952$E(-005)$&3.0775&7.2500$E(-005)$ &3.0372\\
160&4.6748$E(-006)$&3.0212&9.0125$E(-006)$ &3.0080\\
320&5.8184$E(-007)$&3.0062&1.1255$E(-006)$ &3.0014\\
640&7.2629$E(-008)$&3.0020&1.4067$E(-007)$ &3.0001\\
\end{tabular}
\begin{tabular}{rcccc}\hline
 $N$&$\alpha = 5/6$&e.c.r&$\alpha=1$&e.c.r\\[0.5ex] \hline
10&3.3579$E(-001)$&           &3.6830$E(-001)$ &\\
20&9.7882$E(-003)$&5.1004&1.6541$E(-002)$ &4.4768\\
40&9.8787$E(-004)$&3.3087&1.9671$E(-003)$ &3.0719\\
80&1.1104$E(-004)$&3.1533&2.4081$E(-004)$ &3.0301\\
160&1.3330$E(-005)$&3.0583&3.0099$E(-005)$ &3.0001\\
320&1.6404$E(-006)$&3.0226&3.7716$E(-006)$ &2.9946\\
640&2.0374$E(-007)$&3.0092&4.7276$E(-007)$ &2.9960\\
\hline
\end{tabular}
\end{center}
\caption{Errors $\mathrm E_N^{\mathrm{ext}}$ and $\mathrm E_N^\varphi$ for equation \eqref{eq:6.8}, with potential representation \eqref{eq:6.3} (direct, second kind integral equation, Neumann). The upper table corresponds to $\mathrm E_{N}^{\mathrm{ext}}$ and the lower table corresponds to $\mathrm E_{N}^\varphi$.}\label{table:eq7}
\end{table}

\begin{table}[htb]
\begin{center}
\begin{tabular}{rcccc}\hline
 $N$&$\alpha = 5/6$&e.c.r&$\alpha=1$&e.c.r\\[0.5ex] \hline
10&1.7474$E(-001)$&           &1.5968$E(-001)$ &\\
20&7.0648$E(-003)$&4.6284&8.6420$E(-003)$ &4.2077\\
40&4.5988$E(-004)$&3.9413&6.3940$E(-004)$ &3.7566\\
80&4.2497$E(-005)$&3.4358&6.5002$E(-005)$ &3.2982\\
160&4.4521$E(-006)$&3.2548&7.2760$E(-006)$ &3.1593\\
320&5.0449$E(-007)$&3.1416&8.5836$E(-007)$ &3.0835\\
640&5.9863$E(-008)$&3.0751&1.0416$E(-007)$ &3.0428\\
\end{tabular}
\begin{tabular}{rcccc}\hline
 $N$&$\alpha = 5/6$&e.c.r&$\alpha=1$&e.c.r\\[0.5ex] \hline
10&1.0240$E(+000)$&           &9.9018$E(-001)$ &\\
20&1.0248$E(-001)$&3.3207&1.1062$E(-001)$ &3.1621\\
40&5.8839$E(-003)$&4.1225&7.3167$E(-003)$ &3.9182\\
80&4.3947$E(-004)$&3.7429&6.3631$E(-004)$ &3.5234\\
160&3.7752$E(-005)$&3.5411&6.3618$E(-005)$ &3.3222\\
320&3.7134$E(-006)$&3.5457&7.0366$E(-006)$ &2.1765\\
640&4.0517$E(-007)$&3.1962&8.2524$E(-007)$ &3.0920\\
\hline
\end{tabular}
\end{center}
\caption{Errors $\mathrm E_{N}^{\mathrm{ext}}$ and $\mathrm E_{N}^\varphi$ for equation 
\eqref{eq:6.9} with potential representation \eqref{eq:6.3} (direct, hypersingular 
equation, Neumann). The upper table corresponds to $\mathrm E_{N}^{\mathrm{ext}}$ and the 
lower table corresponds to $\mathrm E_{N}^\varphi$.}\label{table:eq8}
\end{table}

\begin{table}[htb]
\begin{center}
\begin{tabular}{rcccc}\hline
 $N$&$\alpha = 5/6$&e.c.r&$\alpha=1$&e.c.r\\[0.5ex] \hline
10&1.2545$E(-001)$&           &1.3462$E(-001)$ &\\
20&1.4698$E(-003)$&6.4153&2.1018$E(-003)$ &6.0011\\
40&2.7686$E(-004)$&2.4084&4.6628$E(-004)$ &2.1723\\
80&4.2480$E(-005)$&2.7043&7.6357$E(-005)$ &2.6104\\
160&5.7877$E(-006)$&2.8757&1.0648$E(-005)$ &2.8422\\
320&7.5027$E(-007)$&2.9475&1.3925$E(-006)$ &2.9348\\
640&9.5326$E(-008)$&2.9765&1.7758$E(-007)$ &2.9711\\
\hline
\end{tabular}
\end{center}
\caption{Errors $\mathrm E_{N}^{\mathrm{ext}}$ for equation \eqref{eq:6.11} with potential representation \eqref{eq:6.10} (indirect, combined field potential, Dirichlet). }\label{table:eq9}
\end{table}

\begin{table}[htb]
\begin{center}
\begin{tabular}{rcccc}\hline
 $N$&$\alpha = 5/6$&e.c.r&$\alpha=1$&e.c.r\\[0.5ex] \hline
10&1.1559$E(-001)$&           &1.2167$E(-001)$ &\\
20&2.7343$E(-003)$&5.4017&3.4492$E(-003)$ &5.1406\\
40&9.5433$E(-005)$&4.8405&1.3958$E(-004)$ &4.6271\\
80&5.6897$E(-006)$&4.0681&9.9258$E(-006)$ &3.8138\\
160&4.0825$E(-007)$&3.8008&6.9688$E(-007)$ &3.8322\\
320&3.2007$E(-008)$&3.6730&5.0311$E(-008)$ &3.7920\\
640&2.9313$E(-009)$&3.4488&4.0423$E(-009)$ &3.6376\\
\hline
\end{tabular}
\end{center}
\caption{Errors $\mathrm E_{N}^{\mathrm{ext}}$ for equation \eqref{eq:6.12} with potential representation \eqref{eq:6.10} (indirect, combined field potential, Neumann). }\label{table:eq10}
\end{table}

\begin{figure}[htb]
\begin{center}
\includegraphics[width=0.45\textwidth]{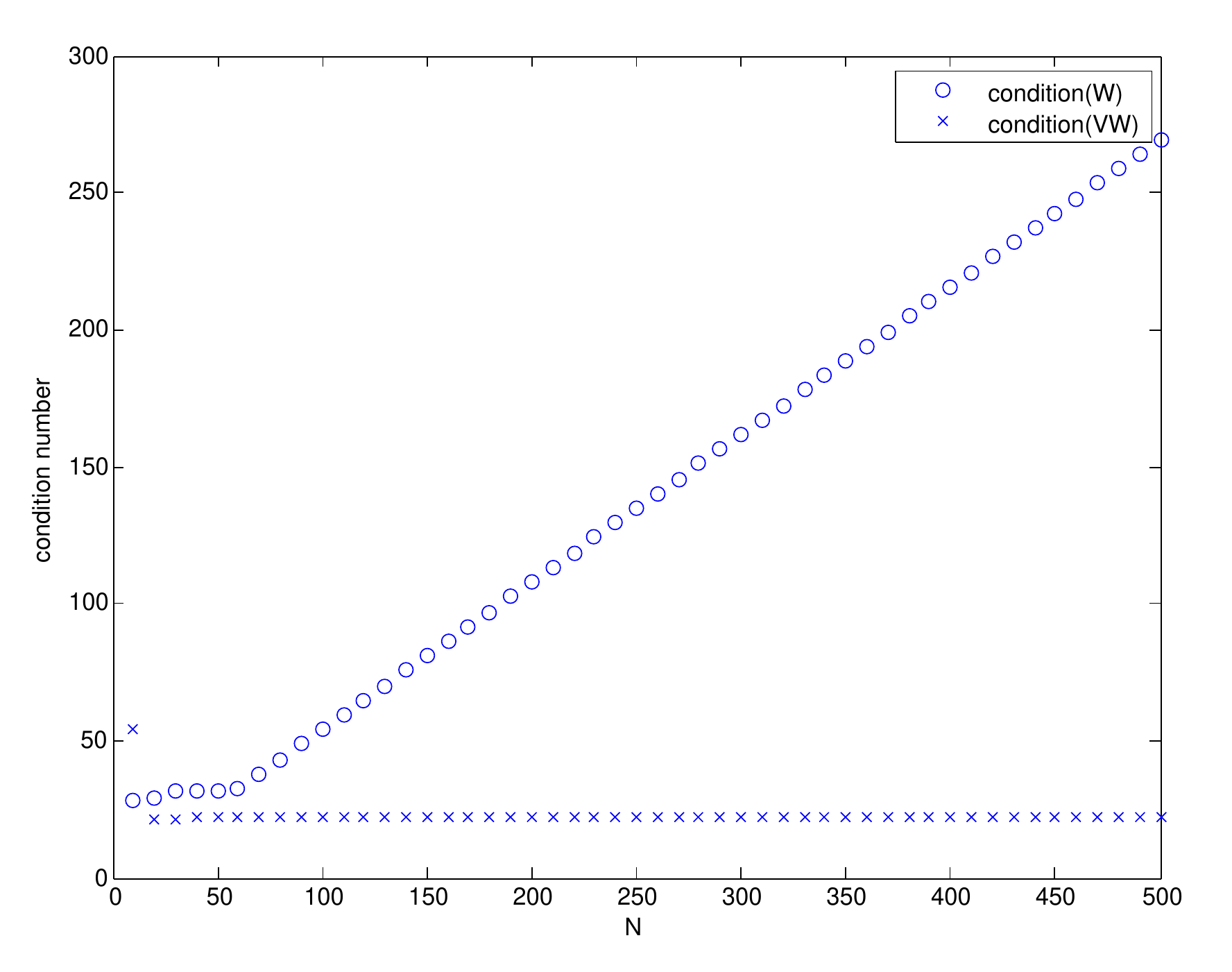}
\includegraphics[width=0.45\textwidth]{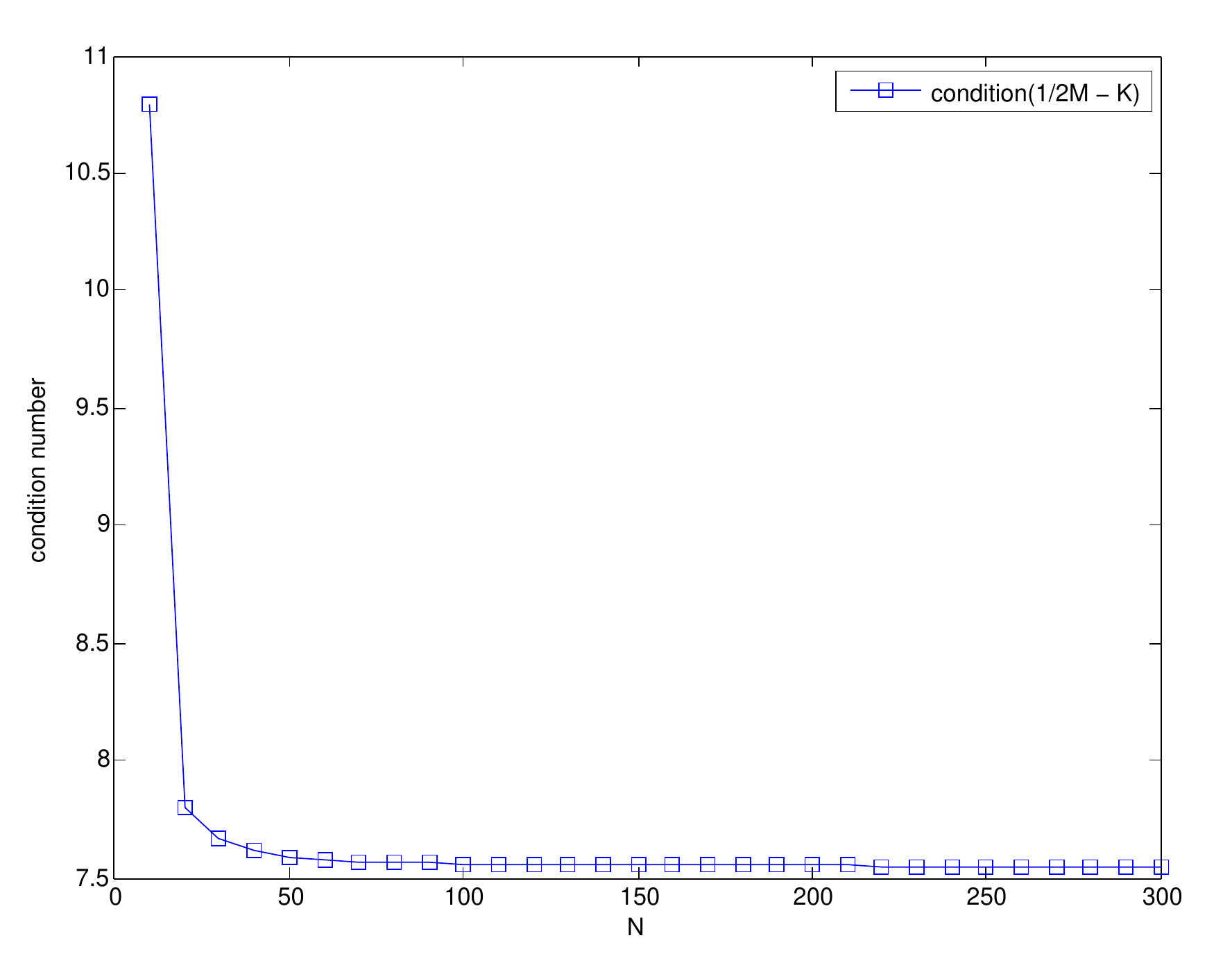}
\end{center}
\caption{Condition numbers for the matrices $\mathrm W_h$, $\mathrm V_h\mathrm W_h$ and $\smallfrac12\mathrm M-\mathrm K_h$. The results are given for the choice $\alpha=5/6$. Results for $\alpha=1$ are almost identical.}
\label{fig:2}
\end{figure}

\paragraph{Tests on condition numbers} Equations associated to weakly singular and 
hypersingular operators will have naturally growing condition numbers. In Figure 
\ref{fig:2} we show how $\mathrm{cond}(\mathrm W_h)=\mathcal O(N)$, but $\mathrm{cond}
(\mathrm V_h\mathrm W_h)=\mathcal O(1)$, that is, the Calder\'on preconditioner works at 
the discrete level. We also show how integral equations of the second kind are well 
conditioned, by showing how $\mathrm{cond}(\smallfrac12\mathrm M-\mathrm K_h)=\mathcal 
O(1).$

\begin{figure}[htb]
\begin{center}
\includegraphics[width=0.45\textwidth]{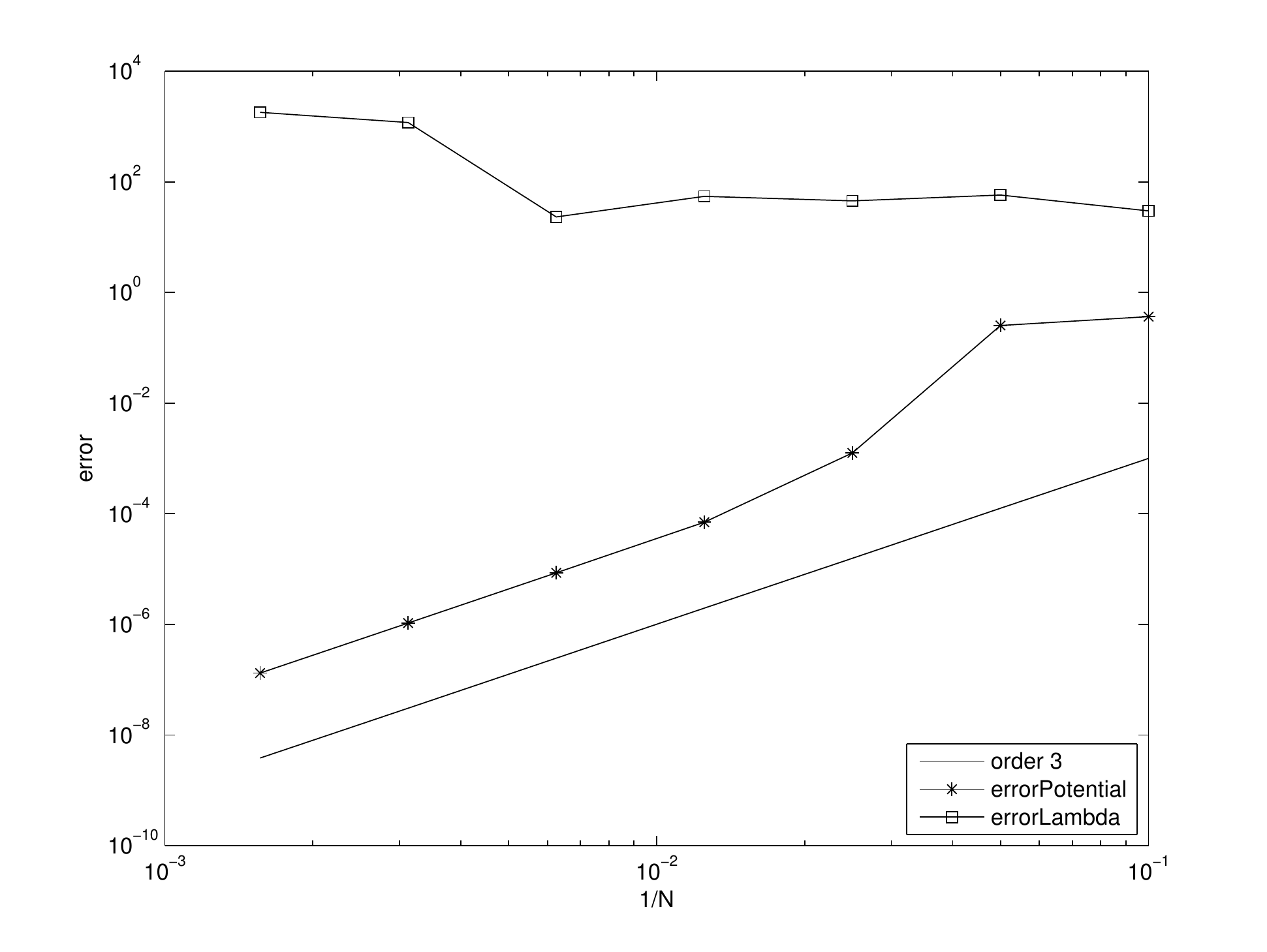}
\includegraphics[width=0.45\textwidth]{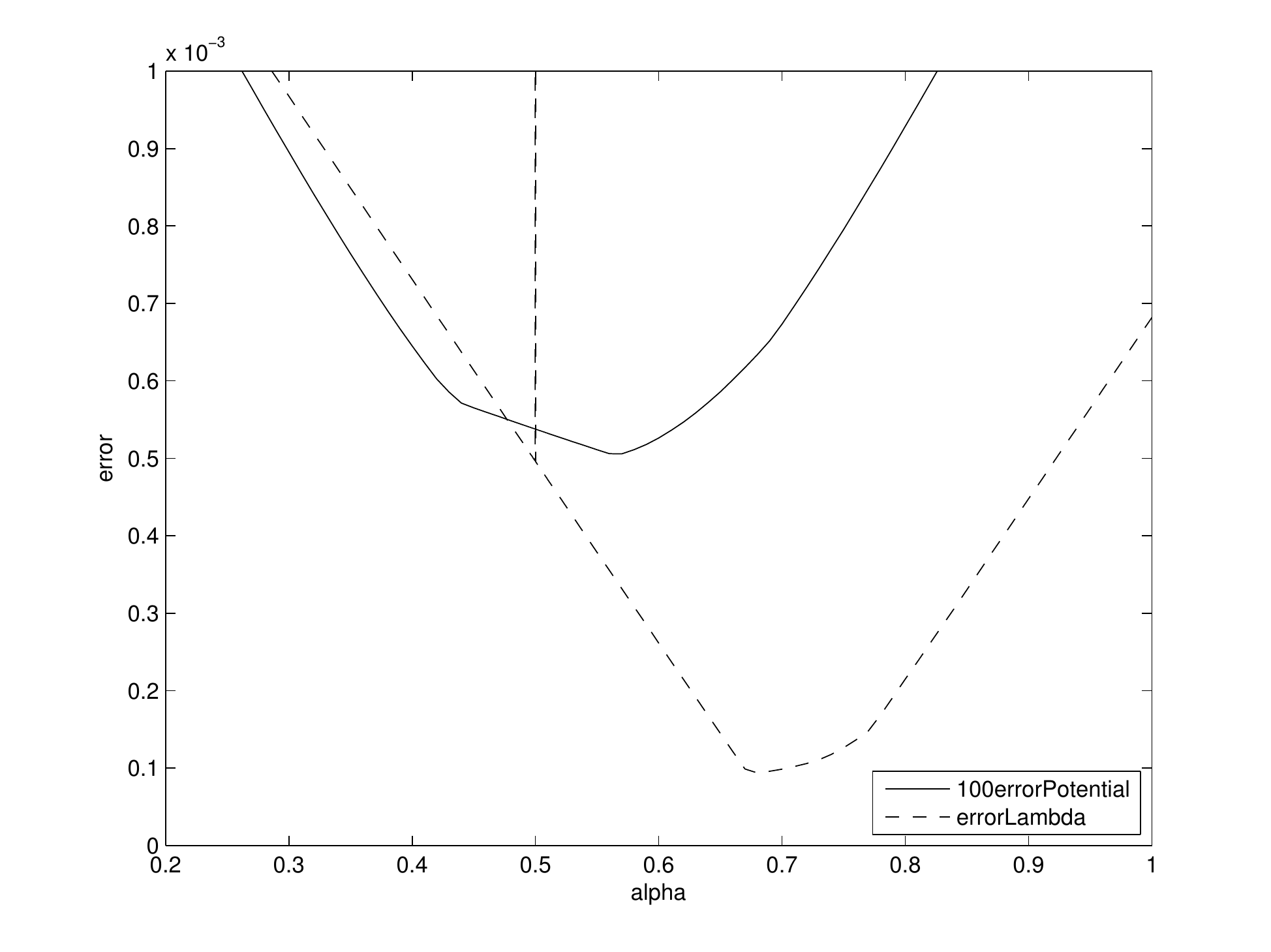}
\end{center}
\caption{The figure on the left shows history of convergence for the choice $\alpha=1/2$ 
using a direct single layer potential based method. The method is clearly not converging 
for the unknown on the boundary, but convergence is restored in the smoothing 
postprocessing of the potential. The figure on the right shows a history of convergence 
w.r.t. $\alpha$ for fixed $N$. The peak at $\alpha=1/2$ corresponds to the unstable choice 
of this parameter.}
\label{fig:alpha}
\end{figure}

\begin{figure}[htb]
\begin{center}
\includegraphics[width=0.45\textwidth]{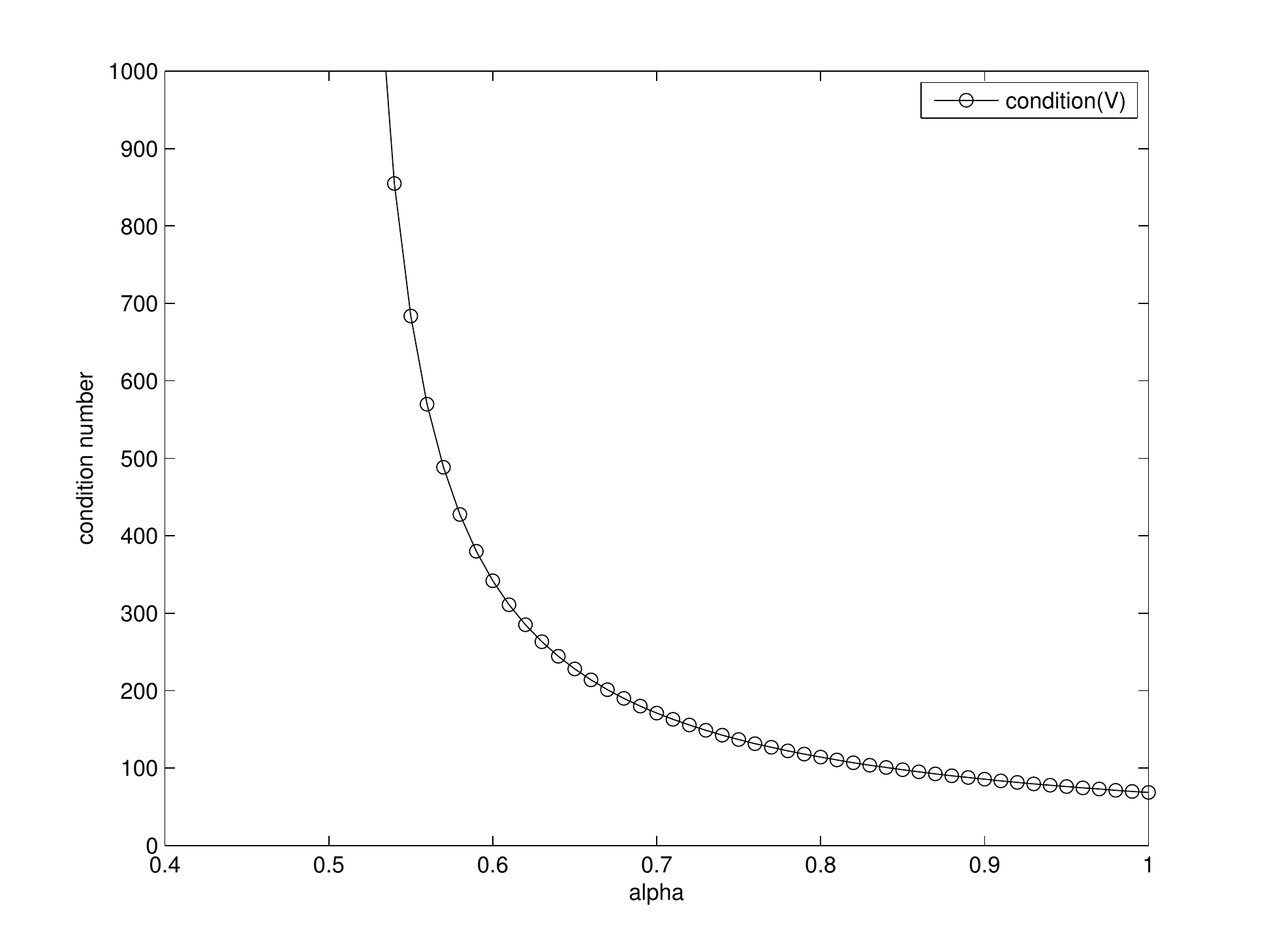}
\includegraphics[width=0.45\textwidth]{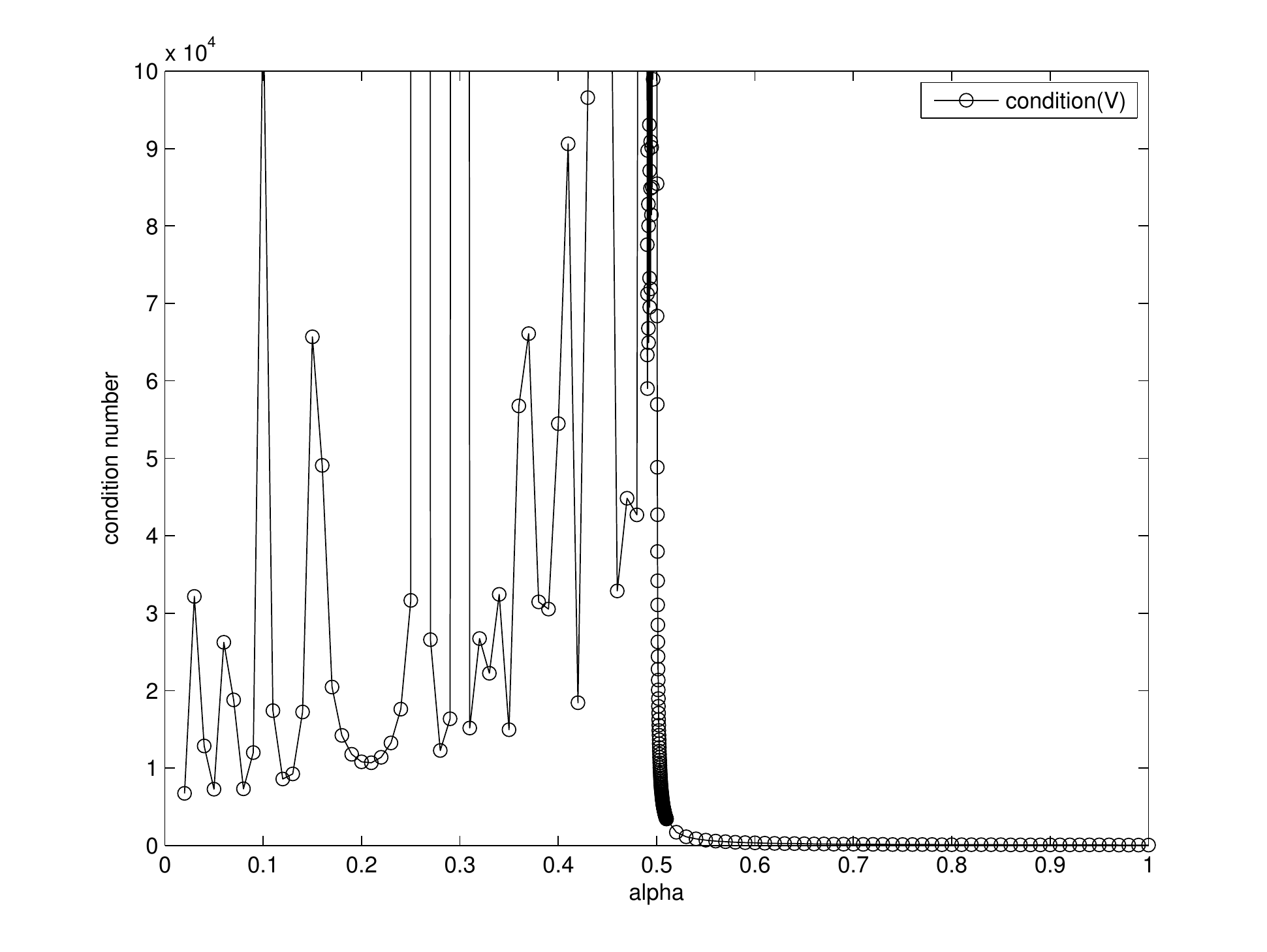}
\end{center}
\caption{Condition number of the matrix $\mathrm V_h$ as a function of the parameter 
$\alpha$. The choice $\alpha=1/2$ equalizes the height of the four Dirac deltas in Figure 
\ref{fig:fork}, making the method unstable. Past this threshold, condition numbers are 
unreasonably high.}
\label{fig:alpha2}
\end{figure}

\paragraph{Dependence with respect to $\alpha$} It is unclear from the experiments whether 
there is a much better choice of the parameter $\alpha$, that dictates the mixture of test 
functions in the method. Let us first show that $\alpha=1/2$ is not feasible. For a test 
equation \eqref{eq:6.4} we compute the errors $\mathrm E_N^\lambda$ and $\mathrm 
E_N^{\mathrm{ext}}$ as $N$ increases. The domain is the curve $\Gamma_1$ and the exact 
solution of the Helmholtz equation is \eqref{eq:7.2}. It is clear from Figure 
\ref{fig:alpha}  that $\mathrm E_N^\lambda$ is not converging, while $\mathrm 
E_N^{\mathrm{ext}}$ converges with the right order. However,  inspection of the condition 
numbers show that they are of the order $10^{20}$. This makes the method highly unstable. 
Convergence of the potential solution can be explained by the fact that the potential 
postprocessing is a smoothing operator which, in some way, eliminates high frequency 
unstable components of the error and only observes approximation properties. In Figure \ref{fig:alpha2}, we explore how the condition numbers of $\mathrm V_h$ blow up as 
$\alpha\to 1/2$ and stay large (but considerably smaller) beyond this value.

\section{More complicated problems}

\subsection{Transmission problems}

\begin{figure}[htb]
\begin{center}
\includegraphics[width=0.45\textwidth]{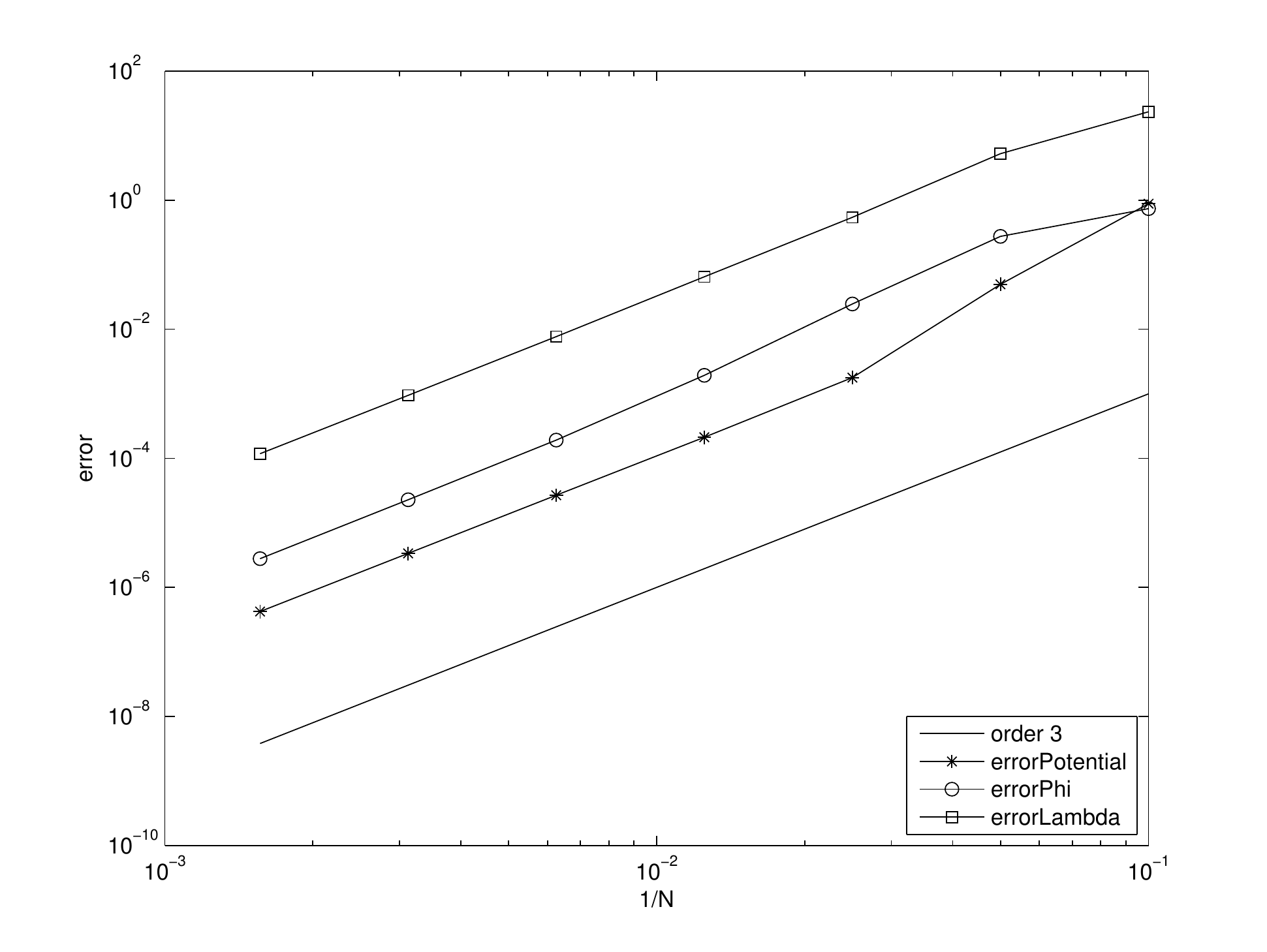}
\includegraphics[width=0.45\textwidth]{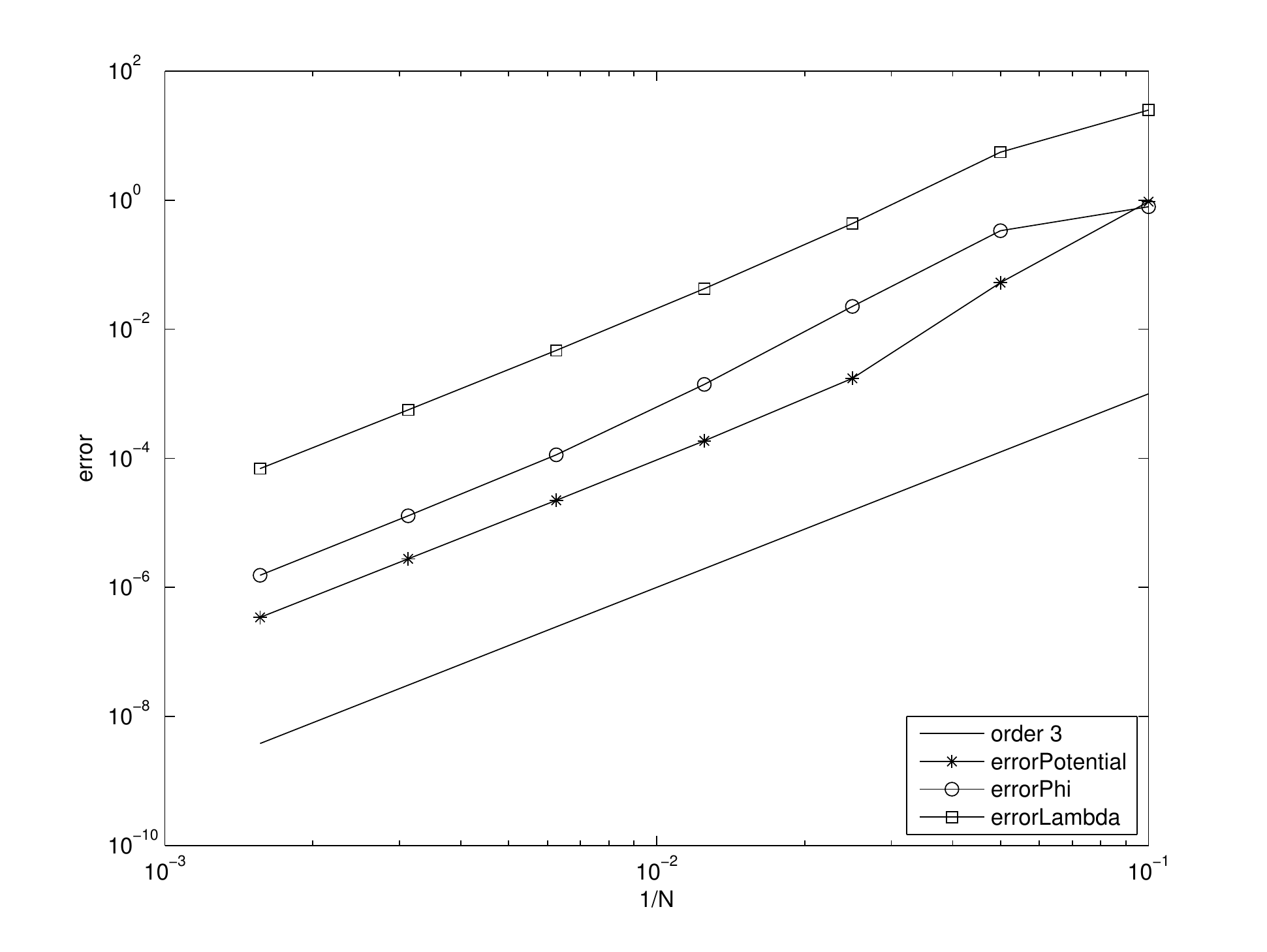}
\end{center}
\caption{Errors for the transmission problem. On the left, errors for the choice $\alpha=1$. On the right, for $\alpha=5/6$.}
\label{fig:1}
\end{figure}

Consider now the domain $\Omega$ interior to the curve \eqref{eq:7.1}. In addition to the exterior Helmholtz equation \eqref{eq:2.2}, we consider an interior equation with a different wave speed
\[
\Delta V+(k/c)^2V=0 \qquad \mbox{in $\Omega$}.
\]
An incident wave $U^{\mathrm{inc}}$ is given and two transmission conditions are imposed on $\Gamma$:
\[
\gamma^+ U+\beta_0
=\gamma^- V, \qquad \partial_{\mathbf n}^+ U+\beta_1=\kappa\,\partial_{\mathbf n}^- V.
\]
In practical problems $(\beta_0,\beta_1):=(\gamma U^{\mathrm{inc}},\partial_{\mathbf n}U^{\mathrm{inc}})$. We choose these transmission data so that the exact solution is the pair given by $U$ in \eqref{eq:7.2} and
\[
V(\mathbf z):=\exp(\imath (k/c)\mathbf z\cdot\mathbf d), \qquad \mathbf d:=(\smallfrac1{\sqrt2},-\smallfrac1{\sqrt2}).
\]
We take $k=3$,  $c=2/3$, and $\kappa=3/2$. The direct symmetric boundary integral formulation of Costabel and Stephan \cite{CoSt:1985} is used. The unknowns are the Cauchy data for the interior problem, so that the integral representations are
\[
U=-\mathrm S(k)(\lambda^--\beta_1)+\mathrm D(k)(\varphi^--\beta_0), \qquad V=\kappa^{-1}\mathrm S(\smallfrac{k}{c}) \lambda^--\mathrm D(\smallfrac{k}{c})\varphi^-.
\]
The corresponding system of integral equations is
\[
\left[\begin{array}{cc} \mathrm W(k)+\kappa\,\mathrm W(\frac{k}{c}) & \mathrm J(k)+\mathrm J(\frac{k}{c}) \\
-\mathrm K(k)-\mathrm K(\frac{k}{c}) & \mathrm V(k)+\kappa^{-1} \mathrm V(\frac{k}{c})\end{array}\right]\left[\begin{array}{c} \varphi^- \\ \lambda^- \end{array}\right]= \left[\begin{array}{cc} \mathrm W(k) & \frac12\mathrm I +\mathrm J(k) \\
\frac12\mathrm I-\mathrm K(k) & \mathrm V(k)\end{array}\right]\left[\begin{array}{c} \beta_0 \\ \beta_1 \end{array}\right].
\]
We discretize each of the elements in the system of integral equations and in the integral representations using the rules of the discrete Calder\'on Calculus. Taking $N$ discretization points on the boundary, we compute the exterior error \eqref{eq:7.3} and errors on the boundary
\[
\mathrm E_N^\lambda:=N\,\max_j|\lambda_j-\kappa\,\nabla V(\mathbf m_j)\cdot\mathbf n_j|, \qquad \mathrm E_N^\varphi:=\max_j| \phi_j-V(\mathbf m_j)|, \quad\mbox{with}\quad\boldsymbol\phi=\mathrm Q\boldsymbol\varphi.
\]
The corresponding errors are plotted in Figure \ref{fig:1}.

\subsection{CQ discretization in the time domain}

In this final example we show how to combine the fully discrete Calder\'on Calculus with a 
Convolution Quadrature routine to produce time-domain discretization of scattering of waves 
by obstacles. We first explain some general ideas of the CQ method. More details, 
specifically applied to scattering problems, are given in \cite{BaSc:2012, LaSa:2009}, while the original ideas of multistep-based CQ for hyperbolic problems appear in  \cite{Lubich:1994}.

\paragraph{Generalities about CQ} We start with a causal approximation of the derivative: if $\kappa>0$, then the operator
\begin{equation}\label{eq:8.1}
\partial_\kappa u:=\smallfrac1\kappa (\smallfrac32 u-u(\cdot-\kappa)+\smallfrac12 u(\cdot-2\kappa))
\end{equation}
is the backward differentiation operator associated to the BDF2 method. The associated transfer function (the Laplace transform of the operator) is
\begin{equation}\label{eq:8.2}
s_\kappa:=\smallfrac1\kappa(\smallfrac32-2e^{-\kappa s}+\smallfrac12 e^{-2\kappa s}).
\end{equation}
Let now $\mathrm A_h(s)$ be any of the elements of the discrete Calculus (one of the potentials or one of the operators), with $k=-\imath s$,  $s\in \mathbb C$ and $\mathrm{Re}\,s>0$. This is the same as saying that we are taking the operators associated to the Laplace resolvent equation $\Delta U-s^2 U=0$ in $\mathbb R^2\setminus\Gamma$ (radiation conditions are reduced to imposing $U\in H^1(\mathbb R^2\setminus\Gamma)$, which in practice imposes exponential decay at infinity). After some manipulation in the complex plane, we can write
\[
\mathrm A_h(s_\kappa)=\sum_{m=0}^\infty \mathrm A_{\kappa,h}[m] e^{-\kappa m s}.
\]
The Convolution Quadrature method is the practical computation of convolutions of the form
\begin{equation}\label{eq:8.3}
\mathrm A_h(\partial_\kappa)\boldsymbol\psi=\sum_{m=0}^\infty \mathrm A_{\kappa,h}[m]\boldsymbol\psi(\cdot-m\kappa)
\end{equation}
(compare with \eqref{eq:8.1} and \eqref{eq:8.2}). The forward convolution form consists of sampling a causal function $\boldsymbol\psi:\mathbb R\to\mathbb C^N$, denoting $\underline{\boldsymbol\psi}[n]:=\boldsymbol\psi(\kappa\,n)$, and then computing
\begin{equation}\label{eq:8.4}
\mathrm A_h(\partial_\kappa) \underline{\boldsymbol\psi}[n]:=\sum_{m=0}^\infty \mathrm A_{\kappa,h}[m]\boldsymbol\psi[n-m].
\end{equation}
(Note that we use the same notation, but now $\underline{\boldsymbol\psi}$ is discrete in time, i.e., it is a sequence of vectors.)
The same idea can be used to solve convolution equations (in the same way that \eqref{eq:8.1} is the seed of the BDF2 method)
\begin{equation}\label{eq:8.5}
\sum_{m=0}^\infty \mathrm A_{\kappa,h}[m]\boldsymbol\psi[n-m]=\boldsymbol\xi[n]\qquad n=0,1,\ldots,
\end{equation}
where $\boldsymbol\xi:\mathbb R\to \mathbb C^N$ is a given causal function sampled at the 
points $\kappa\,n$, or $\boldsymbol\xi[n]$ are the entries of a sequence of vectors 
$\underline{\boldsymbol\xi}$. Note that in \eqref{eq:8.4} and \eqref{eq:8.5} data 
($\boldsymbol\psi$ and $\boldsymbol\xi$ respectively) are sampled in the time domain, while 
the action of the operator is taken using the transfer function. Practical ways of 
computing these convolutions are explained in \cite{BaSc:2012}. They involve a clever use 
of FFT, contour integrals, and multiple evaluations of the transfer function $\mathrm 
A_h(s)$. In the case of the convolution equation \eqref{eq:8.5}, repeated inversion of 
$\mathrm A_{\kappa,h}[0]=\mathrm A_h(s_0)=\mathrm A_h(\smallfrac32 \smallfrac1\kappa)$ is also required.

\paragraph{A scattering problem} In this first example, we use the time domain version of \eqref{eq:6.3} and \eqref{eq:6.9}. The normal derivative of an incident plane wave $U^{\mathrm{inc}}(t,\mathbf x)$  is sampled at the observation points at all times
\[
\boldsymbol\beta_1^\pm[n]:=-(\nabla U^{\mathrm{inc}}(n\,\kappa,\mathbf m_1^\pm)\cdot\mathbf 
n_1,\ldots,\nabla U^{\mathrm{inc}}(n\,\kappa,\mathbf m_N^\pm)\cdot\mathbf n_N)^\top, \quad 
n\ge 0.
\]
We assume that the discrete function $\boldsymbol\beta_1[n]:=\mathrm P^+\boldsymbol\beta_1^+[n]+\mathrm P^-\boldsymbol\beta_1^-[n]$ is causal: this is true in 
the reasonable physical situation when  the incident wave has not reached any of the 
obstacles at time zero. We then solve equations looking for causal sequences $\underline{\boldsymbol\varphi}=(\boldsymbol\varphi[n])$ and $\underline{\boldsymbol\lambda}=(\boldsymbol\lambda[n])$ satisfying
\begin{equation}\label{eq:8.5a}
\mathrm M\boldsymbol\lambda[n]=\boldsymbol\beta_1[n], \qquad \mathrm 
W_h(\partial_k)\underline{\boldsymbol\varphi}[n]=-\smallfrac12\mathrm M\boldsymbol\lambda[n]-\mathrm J_h(\partial_\kappa)\underline{\boldsymbol\lambda}[n], \qquad \forall n\ge 0.
\end{equation}
The potentials are then computed at every time step using the CQ method once again, resulting in  sequences
\begin{equation}\label{eq:8.6}
U[n]=\mathrm S_h(\partial_\kappa)\underline{\boldsymbol\lambda}[n]-\mathrm D_h(\partial_\kappa)\underline{\boldsymbol\varphi}[n].
\end{equation}
Note that this is a fully discrete method for the scattering of a sound-hard obstacle by a 
transient incident wave. Note also that the sequence of functions \eqref{eq:8.6} are a classical solution of the BDF2-discretized wave equation \cite{Lubich:1994}:
\[
\partial_\kappa^2 U[n] -\Delta U[n]=0 \qquad \mbox{in $\mathbb R^2\setminus\Gamma$}\qquad\forall n.
\]
To test the method, we change some signs so that we end up solving an interior boundary 
value problem, namely, we solve $\mathrm 
W_h(\partial_\kappa)\underline{\boldsymbol\varphi}=\frac12\mathrm 
M\underline{\boldsymbol\lambda}-\mathrm 
J_h(\partial_\kappa)\underline{\boldsymbol\lambda}$, instead of the second equation in \eqref{eq:8.5a}. The potential solution \eqref{eq:8.6} is then an approximation of $-U^{\mathrm{inc}}(n\,\kappa,\cdot)$ in $\Omega_-$. 

For the experiments we take the boundary of the domain parametrized with
\[
\smallfrac1{10\sqrt2} (4\,(1+\cos^2(2\pi t))\cos(2\pi t),5\,(1+\sin^2(2\pi t))\sin(2\pi t))\left( \begin{array}{cc} 1 & -1 \\ 1 & 1\end{array}\right),
\]
the incident wave given by
\[
U^{\mathrm{inc}}(t,\mathbf z):=\rho(t\!-\!R\!+\!\mathbf z\cdot\mathbf d), \quad R=1.2, \quad \mathbf d:=(-\smallfrac1{\sqrt2},-\smallfrac1{\sqrt2}), \quad \rho(t):=\sin^3(3\,t)\chi_{t\ge 0},
\]
$N$ discretization points on the curve, and $M$ time steps of length $T/M$, where $T=5$. Finally we compute errors

\begin{eqnarray*}
\mathrm E_{N,M}^{\mathrm{int}} &:=& |U[M](\mathbf z_\circ)+U^{\mathrm{inc}}(T,\mathbf z_\circ)|,\qquad \mathbf z_\circ=(0.2,0.2),\\
\mathrm E_{N,M}^{\varphi} &:=& \max_j |\phi_j[M]+U^{\mathrm{inc}}(T,\mathbf m_j)|, \qquad \boldsymbol\phi[M]=\mathrm Q \boldsymbol\varphi[M].
\end{eqnarray*}
The values of $N$ and $M$ are chosen so that $\mathcal O(N^{-3})=\mathcal O(M^{-2})$: for $j=10,\ldots,19$, we define
\[
N_j =\lfloor 20\, (1.2)^j\rfloor, \qquad M_j:=\lfloor N^{3/2} \,20^{-1/2}\rfloor, \qquad N_j^3 \approx 20 M_j^2, \qquad N_{j+1}/N_j \approx 1.2.
\]
The results are reported in Table \ref{table:10}. Experimental convergence rates are shown to confirm that the errors in $\mathcal O(N^{-3})$.
\begin{table}[htb]
\begin{center}
\begin{tabular}{rccccc}\hline
 $N$& $M$ &$\mathrm E_{N,M}^{\mathrm{ext}}$& e.c.r & $\mathrm E_{N,M}^{\varphi}$&e.c.r\\[0.5ex] \hline
123&305&7.1971$E(-002)$&   &1.3079$E(-001)$&           \\
148&402&4.2559$E(-002)$&2.8816&7.6194$E(-002)$&2.9634\\
178&531&2.4632$E(-002)$&2.9994&4.3811$E(-002)$&3.0353\\
213&695&1.4327$E(-002)$&2.9723&2.5594$E(-002)$&2.9482\\
256&915&8.2648$E(-003)$&3.0173&1.4754$E(-002)$&3.0211\\
308&1208&4.7404$E(-003)$&3.0489&8.4560$E(-003)$&3.0532\\
369&1584&2.7533$E(-003)$&2.9801&4.9135$E(-003)$&2.9776\\
443&2084&1.5894$E(-003)$&3.0135&2.8368$E(-003)$&3.0129\\
532&2743&9.1716$E(-004)$&3.0157&1.6368$E(-003)$&3.0162\\
638&3603&5.3072$E(-005)$&3.0005&9.4818$E(-004)$&2.9946\\
\hline
\end{tabular}
\end{center}
\caption{Errors $\mathrm E_{N,M}^{\mathrm{int}}$ and $\mathrm E_{N,M}^{\varphi}$ for an interior problem in the time domain.}\label{table:10}
\end{table}

\begin{figure}[htb]
\begin{tabular}{cc}
$t=2.8560$  & $t=5.2560$\\
\includegraphics[width=0.45\textwidth]{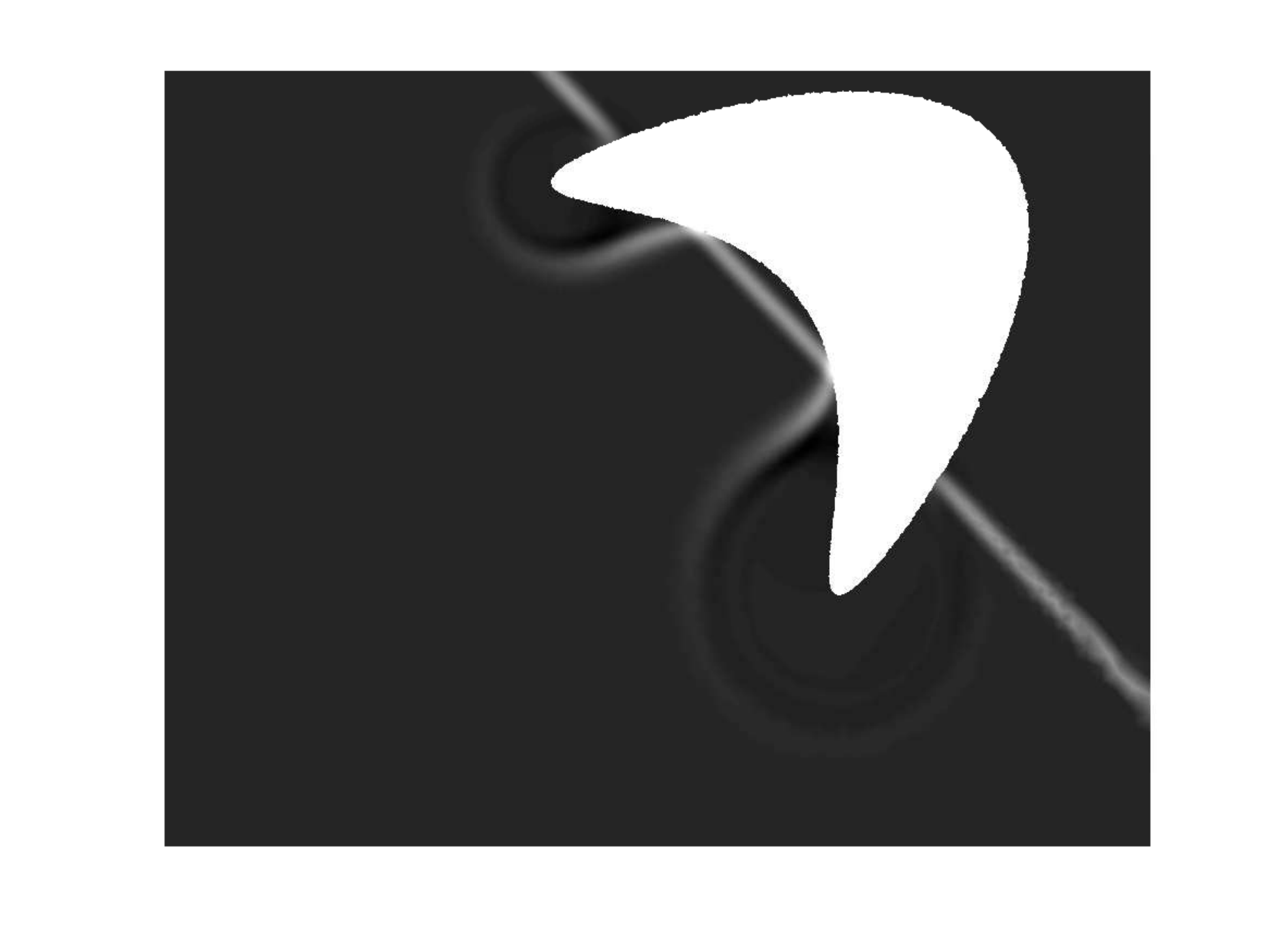} & \includegraphics[width=0.45\textwidth]{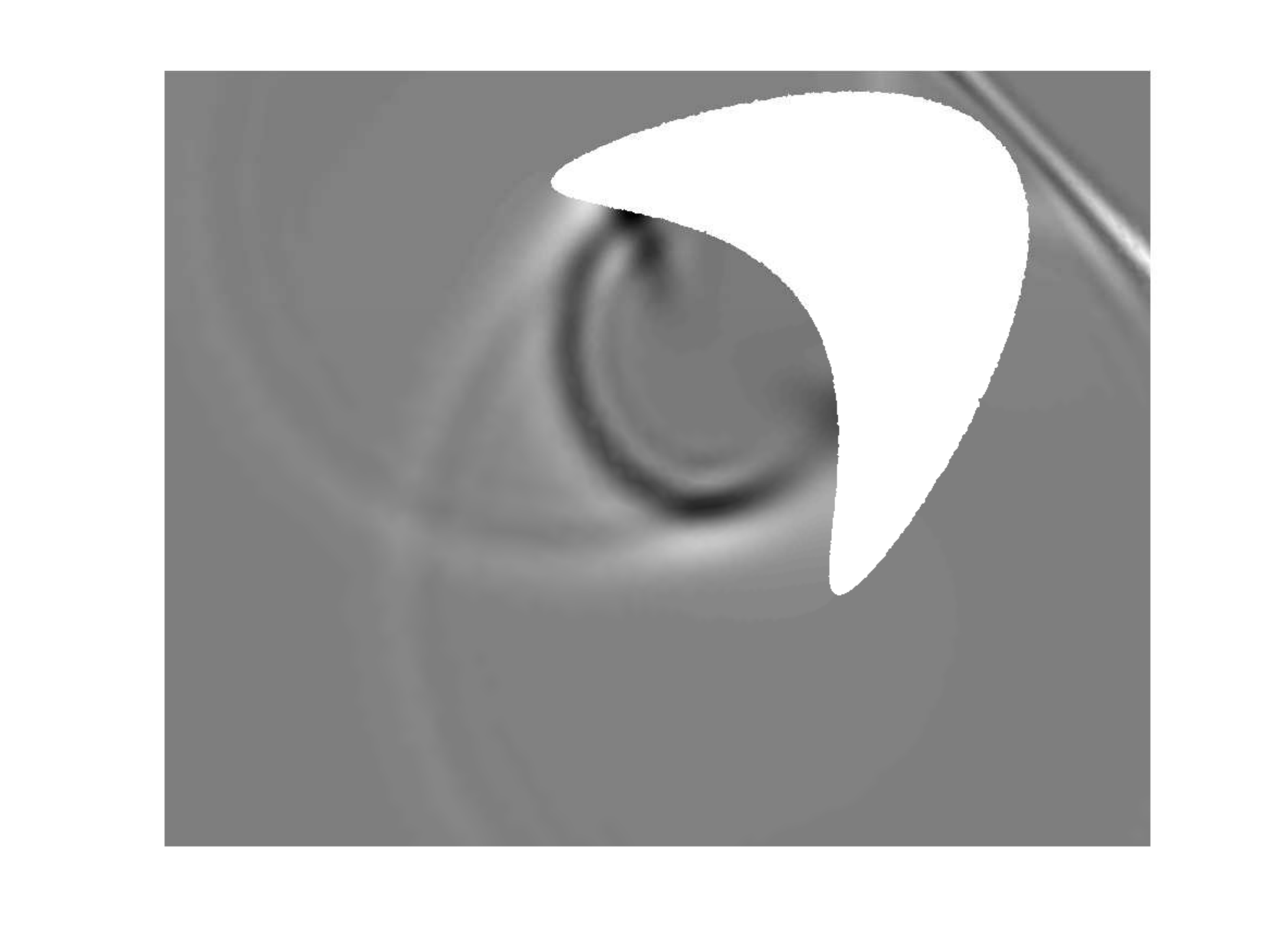} \\
$t=6.6960 $ & $t=7.6560$\\
\includegraphics[width=0.45\textwidth]{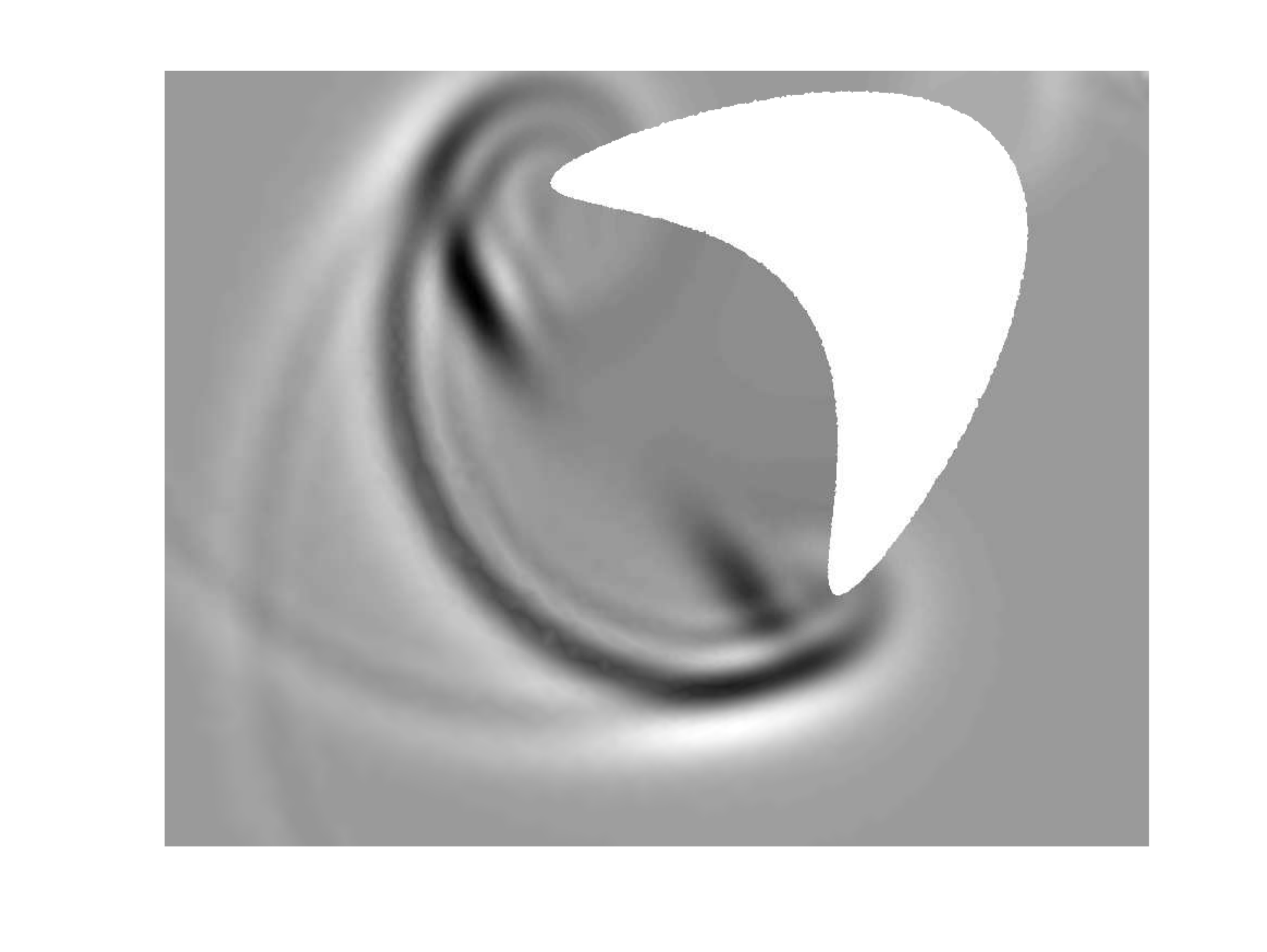} & \includegraphics[width=0.45\textwidth]{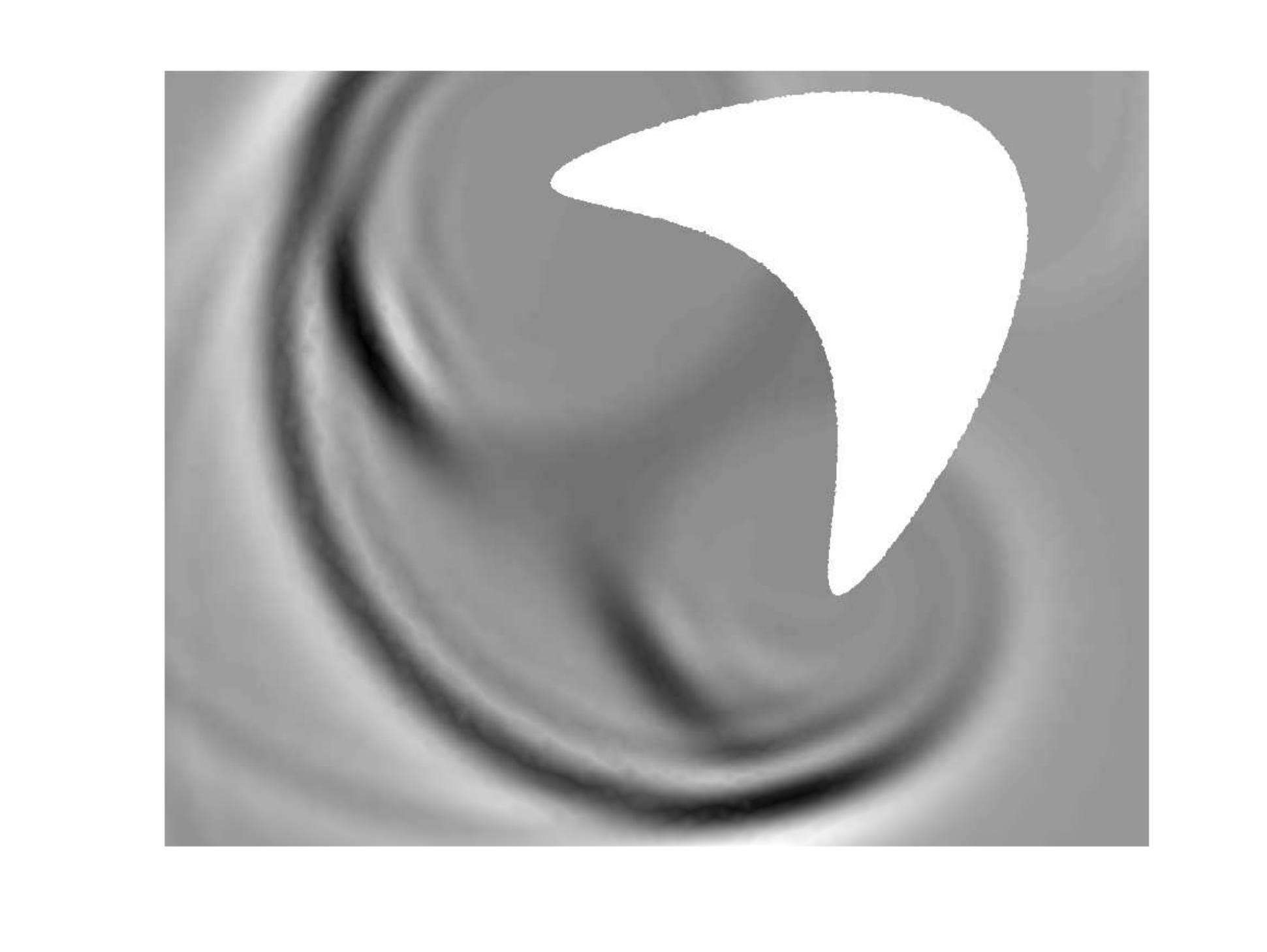}  \\
$t=8.3760$ & $t=9.0960$\\
\includegraphics[width=0.45\textwidth]{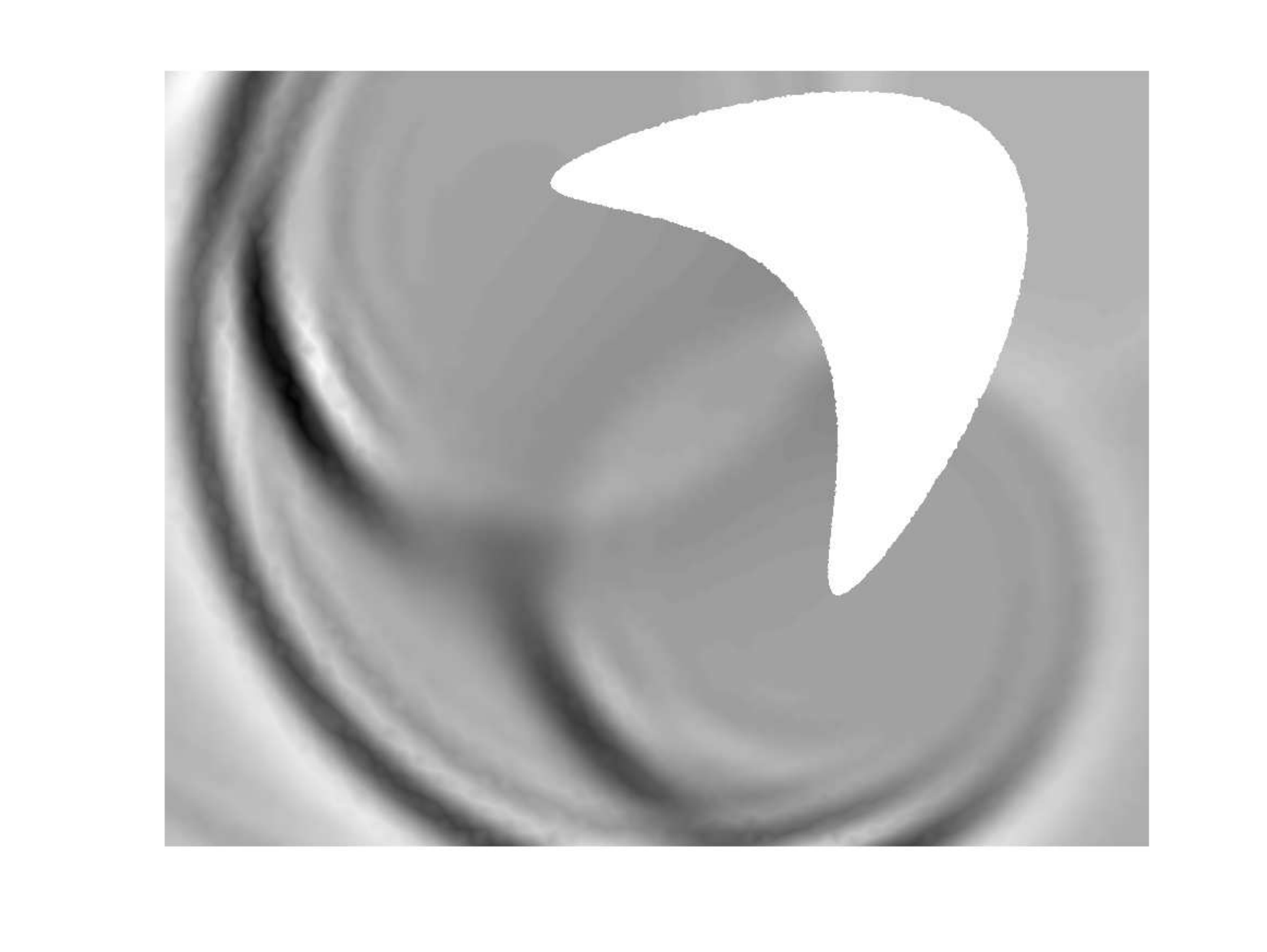} & \includegraphics[width=0.45\textwidth]{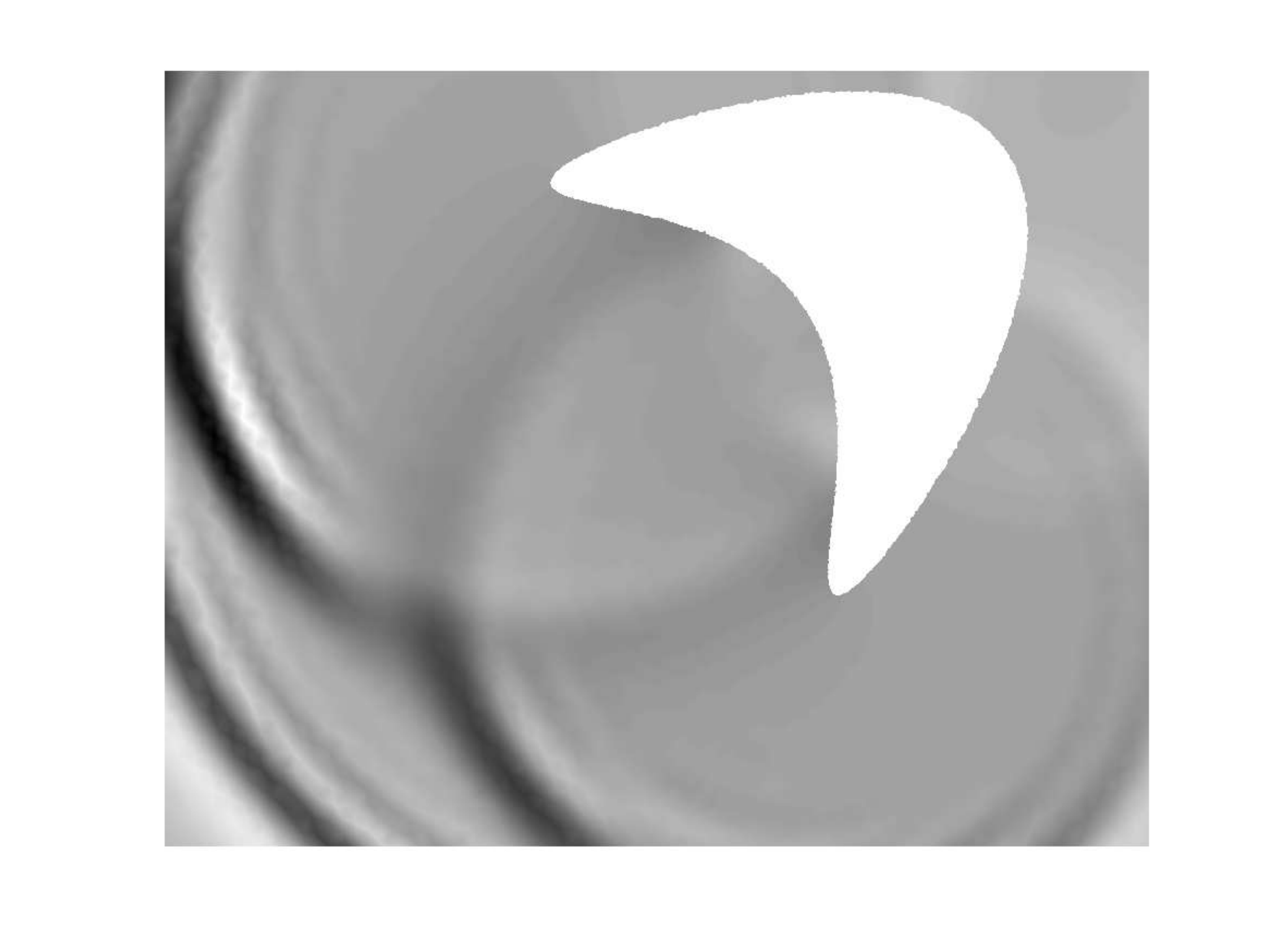}  
\end{tabular}
\caption{Six snapshots of the scattering of a plane wave by a kite-shaped sound-hard 
obstacle. The profile of the wave can be oberserved in the first two images, as it travels 
to the right and to the top. Discretization has been carried out with the order three 
Calder\'on Calculus and a BDF2--based Convolution Quadrature routine.}\label{fig:snapshots}
\end{figure}

\paragraph{A final experiment} To illustrate the capabilities of the time-domain 
discretization, we choose a kite-shaped sound-hard obstacle, hit by a short plane incident 
wave, and we plot several snapshots of the total wave field (incident plus computed wave). Results are shown in Figure \ref{fig:snapshots}.

\bibliographystyle{abbrv}
\bibliography{referencesW}

\end{document}